
\magnification=\magstep1
\def\D{\Delta}

\def\G{\Gamma}

\def\g{\gamma}
\def\ta{\theta}

\def\ds{\displaystyle}
\def\ty{\infty}
\def\l{\left}

\def\r{\right}
\def\o{\over}
\def\ds{\displaystyle}
\def\ss{\scriptstyle}
\def\som#1#2{\sum_{#1}^{#2}}
\def\half{{1\over 2}}

\centerline{{\bf Lecture Notes For An Introductory}}
\centerline{{\bf Minicourse on
 $q$-Series}\footnote{}{{\it 1991 Mathematics Subject
Classification.} Primary 33D15, 33D20,
33D65; Secondary 33D05, 33D45, 33D60, 33D90.}}
\vskip .4cm
\centerline{George Gasper\footnote{}{This work was supported in
part by the National 
Science Foundation under grant DMS-9401452.}}
\centerline{ (September 19, 1995 version)}

\vskip .4cm
These lecture notes were written for a mini-course that was
designed
to introduce students
and researchers to  {\it $q$-series,} 
which are also called {\it basic hypergeometric series} because
of the
parameter $q$ that is used as a base in series that are ``{\it
over, above
or beyond}'' the {\it geometric series}.
We start by considering  $q$-extensions (also 
called $q$-analogues) of the 
binomial theorem, the exponential and gamma functions, and of the
beta 
function and beta integral, and then progress on to the
derivations
of rather general summation, transformation, and expansion
formulas, 
 integral representations, and applications. Our main emphasis is
on 
 methods that can be used to
{\bf derive} formulas, rather than 
to just {\it verify} previously derived formulas.  Since the best
way to
{\it learn} mathematics is to {\it do} mathematics,
in order to enhance the
learning process and enable the reader to practice with the
discussed
methods and formulas, we have provided several carefully selected
exercises at the end
of each section.
We strongly encourage you to obtain a deeper
 understanding of $q$-series by looking at these exercises
and doing {\it at least three} of them in each section. 
Solutions to several of the exercises are given in 
the Gasper and Rahman [1990a] 
``Basic Hypergeometric Series'' book (which we will refer to as
BHS) 
along with additional exercises
 and material on this subject, and references to applications 
to affine root systems (Macdonald identities),  Lie algebras and
groups,
 number theory, orthogonal polynomials, physics (such as
representations of quantum groups and Baxter's work on the
hard hexagon model), statistics,
etc. In particular, for applications to
orthogonal polynomials we recommend the Askey and Wilson [1985]
A.M.S. Memoirs.
For applications to number theory, physics and related fields,
 we recommend the Andrews [1986]
and Berndt [1993] lecture notes and the Fine [1988] book.  To add
a historical
perspective, we have followed the method in BHS of referring to
papers
and books in the References at the end of these notes by placing
the year of publication in square
brackets immediately after the author's name.
\bigskip
\centerline{\bf  1 The $q$-binomial theorem and related formulas}
\bigskip
{\bf  1.1 The  binomial theorem.}  One of the most 
elementary summation formulas
for power series is the sum of the geometric series
$$\sum_{n=0}^\infty z^n = (1 - z)^{-1}, \quad\eqno (1.1.1)$$
where $z$ is a real or complex number and  $|z| < 1.$ 
 By  applying Taylor's theorem to the function $f(z)=(1 -
z)^{-a},$
which is an analytic function of $z$ for $|z| < 1,$ 
and observing by mathematical induction that 
$f^{(n)}(z)=(a)_n (1 - z)^{-a-n}|_{z=0} = (a)_n,$ where $(a)_n$
is the 
{\it shifted factorial} defined by
$$(a)_0 =1,\; (a)_n = a(a+1) \cdots (a+n-1) = {\Gamma
(a+n)\over\Gamma (a)},\quad n=1,2,\ldots, \eqno (1.1.2)$$
one can extend (1.1.1) to 
$$\som{n=0}{\infty}
{(a)_n\over n!}z^n = (1-z)^{-a}, \qquad|z| <1. \eqno (1.1.3)$$
This formula is usually called the {\it binomial theorem}
because, when
$a=-m$ is a negative integer and $z=-x/y,$ it reduces to the
binomial theorem
for the $m$-th power of the binomial $x+y$:
$$(x+y)^m =\som{n=0}{m}
{m\choose n} x^n y^{m-n},\quad m=0,1,2,\ldots \ . \eqno (1.1.4)$$
\bigskip
{\bf  1.2 The $q$-binomial theorem.} Let $0<q<1.$ Since, by
l'H\^opital's rule, 
$$\lim_{q\rightarrow 1^-} {1-q^a\over 1-q} = a \eqno (1.2.1)$$
and hence
$$ \lim_{q\rightarrow 1^-} {(1-q^a)(1-q^{a+1})\cdots
(1-q^{a+n-1})\over (1
-q) (1-q^2)\cdots
(1-q^n)} = {(a)_n\over n!},\;\eqno (1.2.2)$$
it is natural to consider what happens when the
 coefficient of each $z^n$ in (1.1.3) is replaced by the ratio
displayed
on the left side of (1.2.2)
or, more generally, by
$${(a;q)_n\over(q;q)_n},$$
where $(a;q)_n$ is the  $q$-{\it shifted factorial} defined by
$$(a;q)_n = \cases{1,&$n=0,$\cr(1-a)(1-aq)\cdots
(1-aq^{n-1}),&$n=1,2,\ldots$ .\cr}\eqno (1.2.3).$$
\par
Hence, let us set
$$f(a, z)=
 \som{n=0}{\infty} {(a;q)_n\over(q;q)_n}z^n,\eqno (1.2.4)$$
where $a$, $q$, $z$ are real or complex numbers such that $|z|<1$
and, unless
stated otherwise, it is assumed that $|q|<1.$  The case when
$|q|>1$ will
be considered later.  Note that, by the ratio test,  since
$|q|<1$
the  series in  (1.2.4) converges for $|z|<1$ to a function,
which we have
denoted by $f(a, z).$  One 
way to find a formula for $f(a, z)$ analogous to that for the sum
of the 
 series in (1.1.3) is to first observe that,
 since $1-a = (1-a q^n)-a(1-q^n),$
$$\eqalignno{f(a, z)&= 1+\som{n=1}{\infty}
{(a;q)_n\over(q;q)_n}z^n\cr&=
1+\som{n=1}{\infty} {(a q;q)_{n-1}\over(q;q)_n}[(1-a
q^n)-a(1-q^n)]z^n\cr&=
1+\som{n=1}{\infty} {(a q;q)_{n}\over(q;q)_n}z^n-
a\som{n=1}{\infty}{(a q;q)_{n-1}\over(q;q)_{n-1}}z^n\cr&=  f(a q,
z)-
a z f(a q, z) =(1-a z) f(a q, z). &(1.2.5)\cr}$$
By iterating this functional equation $n-1$ times we get that
$$f(a, z)=(a z; q)_n \ f(a q^n, z),\qquad n=1,2,\ldots ,$$
which on letting $n\rightarrow \infty$ and using $q^n \rightarrow
0$ gives
$$f(a, z)=(a z; q)_\infty \ f(0, z) \eqno (1.2.6)$$
with $(a ; q)_\infty $ defined by
$$(a;q)_\infty = \lim_{n\to\ty}(a ; q)_n =\prod^\ty_{k=0} 
(1-aq^k), \qquad  |q| <1. \eqno (1.2.7)$$
Since the above infinite product diverges when $a\ne 0$ and $|q|
\ge 1$, 
whenever $(a;q)_\ty$ appears in a formula, it is usually assumed
that $|q| < 1$.
Now set $a=q$ in  (1.2.6) to obtain 
$$f(0, z)={f(q, z)\o (q z; q)_\infty} = {(1 - z)^{-1} \o (q z;
q)_\infty} =
{1 \o ( z; q)_\infty},$$
which, combined with (1.2.6) and (1.2.4), shows that
$$\som{n=0}{\infty} {(a;q)_n\over(q;q)_n}
z^n = {(az;q)_\infty\over (z;q)_\infty},\quad |z|<1,\ |q|<1.\eqno
(1.2.8)$$
This summation formula was derived by Cauchy [1843] and Heine
[1847].  
It is called the $q$-{\it binomial theorem}, because it is a 
$q$-analogue of the binomial theorem in the sense that
$$\lim_{q\rightarrow 1^-} \ \som{n=0}{\infty}
{(q^a;q)_n\over(q;q)_n}z^n 
=\som{n=0}{\infty}
{(a)_n\over n!}z^n = (1-z)^{-a}, \qquad |z| <1. \eqno (1.2.9)$$
\par
Notice that (1.2.8) and (1.2.9) yield
$$\lim_{q\to 1^-} {(q^az;q)_\ty\o (z;q)_\ty} = (1-z)^{-a},\quad
|z|< 1,\quad a\;{\rm real},\eqno (1.2.10)$$
which, by analytic continuation, holds for $z$ in the complex
plane
cut along the positive real axis from 1 to $\ty$, with $(1-z)^{-
a}$ positive when $z$ is real and less than 1.
In the special case
 $a=q^{-m}, m=0,1,2,\ldots ,$ (1.2.8) gives
$$\som{n=0}{m} {(q^{-m};q)_n\over(q;q)_n}z^n = (zq^{-m};q)_m,
\eqno (1.2.11)$$
where, by analytic continuation, $z$ can be any complex
number.
\par
Heine's proof of the $q$-binomial theorem, which is presented in
Heine [1878],
Bailey [1935, p. 66], Slater [1966, p. 92], and in
 \S 1.3 of BHS 
along with some motivation from Askey [1980], consists of using
series manipulations to derive two difference equations that are
then used
to derive the functional equation
$$ (1-z) f(a,z) = (1-az) f(a,qz). \eqno (1.2.12)$$
Iterating this equation $n-1$ times and letting $n\rightarrow
\infty$ 
gives
$$\eqalignno{f(a,z)  &= {(az;q)_n\over (z;q)_n} f(a,q^nz)\cr&= 
{(az;q)_\infty\over
(z;q)_\infty} f(a,0) = {(az;q)_\infty\over(z;q)_\infty}
&(1.2.13)\cr}$$
and completes the proof.  Since it is easy to use series
manipulations
to show that the  series in (1.2.4) satisfies the functional
equation
(1.2.12), once (1.2.12) has been
{\it discovered} it can be used to give a short {\it
verification} type
proof of the $q$-binomial theorem.
\par
Another derivation of the $q$-binomial theorem can be given by
 calculating the coefficients $c_n=g_a^{(n)}(0)/n!,
n=0,1,2,\ldots,$
in the Taylor series expansion of the
function
$$g_a(z) ={(az;q)_\infty\over (z;q)_\infty}=\som{n=0}{\infty}c_n
z^n
,\eqno(1.2.14)$$
which is an analytic 
function of $z$ when $|z|<1$ and $ |q|<1.$ Clearly $c_0
=g_a(0)=1.$ 
One may show that $c_1 =g_a^{\prime}(0)=(1-a)/(1-q)$ 
by taking the logarithmic derivative of ${(az;q)_\infty\over
(z;q)_\infty}$
and then setting $z=0.$  But, unfortunately, the succeeding
higher order
derivatives of $g_a(z)$ become more and more difficult to
calculate,
and so one is forced to abandon this approach and to search for 
another way to calculate
all of the $c_n$ coefficients.  One simple
method is to notice that from the definition of $g_a(z)$ as
the quotient of two infinite products, it immediately
follows that  $g_a(z)$ satisfies the functional equation
$$ (1-z) \; g_a(z) = (1-az) \; g_a(qz), \eqno (1.2.15)$$
which is of course the same as the functional equation (1.2.12)
 satisfied by $f(a,z).$ In a {\it verification} type proof of
the $q$-binomial theorem, (1.2.15) provides
{\it substantial} motivation for showing, as in Heine's proof,
that the sum of the $q$-binomial series $f(a,z)$
satisfies the functional equation (1.2.12).
\par
To calculate the  $c_n$ coefficients, we first use (1.2.15) to
find that
$$\som{n=0}{\infty}c_n \, z^n -\som{n=0}{\infty}c_n \, z^{n+1} =
\som{n=0}{\infty}c_n \, q^n \ z^n-a\som{n=0}{\infty}c_n \, q^n \,
z^{n+1},$$
or, equivalently,
 $$1+\som{n=1}{\infty}(c_n-c_{n-1}) \, z^n  =1+
\som{n=1}{\infty}(c_n q^n -a c_{n-1}\,  q^{n-1})\,  z^{n},$$
which implies that
$$c_n-c_{n-1} =c_n q^n -a c_{n-1}\, q^{n-1}$$
and hence
$$c_n= {{1-aq^{n-1}}\o{1-q^n}} \, c_{n-1} \eqno (1.2.16)$$
for $n=1,2,\ldots \ .$ \ \
Iterating the recurrence relation (1.2.16) \ $n-1$ times gives
$$c_n= {(a;q)_n\over(q;q)_n}  c_0 ={(a;q)_n\over(q;q)_n}, 
\qquad n=0, 1,2,\ldots, \eqno (1.2.17)$$
which concludes our third derivation of the $q$-binomial theorem
(1.2.8).
For a combinatorial proof using a bijection between two classes
of partitions,
see Andrews [1969].
\bigskip
{\bf  1.3 Related formulas.} 
One immediate consequence the $q$-binomial theorem is the product
formula
$$\som{j=0}{\infty} {(a;q)_j\over(q;q)_j}\, z^j \
\som{k=0}{\infty} {(b;q)_k\over(q;q)_k}\ (az)^k =
\som{n=0}{\infty} {(ab;q)_n\over(q;q)_n}\, z^n,  
\qquad|z| <1,\ |q| <1,\eqno (1.3.1)$$
 which is a $q$-analogue of
 $(1-z)^{-a} (1-z)^{-b} =
(1-z)^{-a-b}.$  By setting $j=n-k$ in the product on the left
side of
(1.3.1) and comparing the coefficients of $z^n$ on both sides of
the
equation, we get
$${(ab;q)_n\over(q;q)_n}=\som{k=0}{n} {(a;q)_{n-k}(b;q)_k
\over(q;q)_{n-k}(q;q)_k}  \, a^k, \eqno (1.3.2)$$
which gives a $q$-analogue of (1.1.4) in the form
$$(ab;q)_n=\som{k=0}{n} \l[{n\atop k}\r]_q (a;q)_{n-k}(b;q)_k \,
a^k
, \eqno (1.3.3)$$
where the $q$-{\it binomial coefficient} is defined
by
$$\l[{n\atop k}\r]_q = {(q;q)_n\o (q;q)_k(q;q)_{n-k}},
\qquad k = 0,1, \ldots, n. \eqno (1.3.4)$$
\par
If we let
$$e_q(z) = \som{n=0}{\infty} {z^n\over
(q;q)_n}, \qquad\;|z| <1,\eqno (1.3.5)$$
then the case $a=0$ of (1.2.8) gives
$$e_q(z) = {1\over (z;q)_\infty} , \qquad\;|z| <1,\ |q|<1.\eqno
(1.3.6)$$
The function $e_q(z)$ is a $q$-analogue of the exponential
function $e^z,$ 
since 
$$\lim_{q\to1^-} e_q(z(1-q)) =e^z.\eqno (1.3.7)$$
Another  $q$-analogue of $e^z$ can be obtained from (1.2.8) by
replacing $z$ with $-z/a$ and then letting $a\rightarrow \infty$
to
find that the $q$-{\it exponential function} defined by
$$E_q(z) =\som{n=0}{\infty}
{q^{n(n-1)/2}\over(q;q)_n}\,  z^n, \qquad\;|z| <\ty, \
|q|<1,\eqno (1.3.8)$$
equals $(-z;q)_\infty$ for all complex values of $z$ and
satisfies
the limit relation
$$ \lim_{q\to1^-} E_q(z(1-q))
=e^z, \qquad |z|<\ty.\eqno (1.3.9)$$
Hence, $e_q(z)E_q(-z) = 1$ when $ |z| <1$ and $|q|<1.$
\par
In Exercise 1.1 below you will be asked to verify the {\it
inversion identity}
$$( a;q)_n = (a^{-1};q^{-1})_n (-a)^nq^{{\scriptstyle
n\mathstrut\choose \scriptstyle 2}}\eqno (1.3.10)$$
for $n=0,1,2,\ldots.$   This 
identity enables us to convert a $q$-series 
formula containing sums of quotients of products
of $q$-shifted factorials in base $q$ to a similar formula with
 base $q^{-1}.$  In particular, it follows from
(1.2.8) and (1.3.10) that if
$|q|>1,$ then the $q$-binomial theorem takes the form
$$\som{n=0}{\infty} {(a;q)_n\over(q;q)_n}\ z^n = 
{(z/q;q^{-1})_\infty\over (az/q;q^{-1})_\infty},\quad |az/q|<1,\
|q|>1.\eqno
(1.3.11)$$
\bigskip
\centerline{\bf Exercises 1}
\bigskip
\item{1.1} Verify the inversion identity (1.3.10) and show
for nonnegative integers $k$  and $n$ that
$$\leqalignno{(a;q)_n &= {(a;q)_\infty\over (aq
^n;q)_\infty},&(i)\cr
(a;q)_{n+k} &= (a;q)_n (aq^n;q)_k,&(ii)\cr
(aq^n;q)_k &= {(a;q)_k (aq^k;q)_n\over (a;q)_n},&(iii)\cr
(aq^{-n};q)_k &= {(a;q)_k (qa^{-1};q)_n\over
(a^{-1}q^{1-k};q)_n}\,  q^{-nk},&(iv)\cr
{1-aq^{2n}\over 1-a} &=
{(qa^\half;q)_n (-qa^\half;q)_n\over (a^\half;q)_n
(-a^\half;q)_n}.&(v)\cr}$$
\medskip
\item{1.2} Use series manipulations to verify that the function
$f(a,z)$
defined in (1.2.4)
satisfies the functional equation (1.2.12).
\medskip
\item{1.3}  Show that if we define  $e_q(z)$ for $|q|>1$ by
$$e_q(z) = \som{n=0}{\infty} {z^n\over
(q;q)_n}, \qquad\;|z| <\ty,\ |q|>1,$$
 then it follows from (1.3.8) and (1.3.10) that
$e_{q^{-1}}(z)=E_q(-qz)$ when  $|z| <\ty$ and $|q|<1. $
\medskip
\item{1.4}  Derive (1.3.11) by using one (or more) of the 
methods used in
\S 1.2 to derive (1.2.8).
\medskip 
\item{1.5} Investigate what happens if, as in the 
third derivation of the $q$-binomial theorem presented in \S 1.2,
you try to calculate the $c_n$ coefficients in (1.2.14) by
starting with the
functional equation $ g_a(z) = (1-az)\; g_{aq}(z),$ which
corresponds to 
the functional equation derived in (1.2.5), 
instead of with the more complicated
functional equation (1.2.15).
\medskip
\item{1.6}  Extend the definition of the  
$q$-binomial coefficient (1.3.4) to
$$\l[{\alpha\atop\beta}\r]_q = {(q^{\beta+1};q)_\ty
(q^{\alpha-\beta +
1};q)_\ty\o (q;q)_\ty (q^{\alpha + 1};q)_\ty}$$
\item{}
for complex $\alpha$ and $\beta$ when $|q| < 1$.  Show that
$$\leqalignno{
\l[{\alpha\atop k}\r]_q &= {(q^{-\alpha};q)_k\o (q;q)_k} \l(-
q^\alpha\r)^kq^{-{\scriptstyle k\mathstrut\choose\scriptstyle
2}},&(i)\cr
\l[{k+\alpha\atop k}\r]_q &= {(q^{\alpha+1};q)_k\o
(q;q)_k},&(ii)\cr
\l[{\alpha+1\atop k}\r]_q &= \l[{\alpha\atop k}\r]_q q^k +
\l[{\alpha\atop k-1}\r]_q = \l[{\alpha\atop k}\r]_q +
\l[{\alpha\atop k-1}\r]_q q^{\alpha + 1-k},&(iii)\cr}$$
\item{} where $k$ and $n$ are nonnegative integers.
\medskip
\item{1.7} Use Ex.$\,$1.6 ($iii$) and induction to prove that if
$x$ and $y$ are
indeterminates such that $xy = qyx, q$ commutes with $x$ and
$y$, and the associative law holds, then
$$(x+y)^n = \sum^n_{k=0} \l[{n\atop k}\r]_q y^k x^{n-k} =
\sum^n_{k=0} \l[{n\atop k}\r]_{q^{-1}} x^k y^{n-k}.$$
See, e.g., Koornwinder[1989].
\bigskip
\centerline{\bf  2 Ramanujan's $_1\psi_1$ summation formula,
Jacobi's triple}
\centerline{\bf product identity and theta functions}
\bigskip
{\bf  2.1  Ramanujan's $_1\psi_1$ summation formula.}   Since
$$(a;q)_n = {(a;q)_\infty\over (aq
^n;q)_\infty}, \qquad |q|<1, \ n=0,1,2,\ldots,\eqno (2.1.1)$$
we may extend the definition (1.2.3) of $(a;q)_n$ to
$$(a;q)_\alpha = {(a;q)_\infty\over (aq
^\alpha;q)_\infty},\qquad |q|<1, \quad\eqno (2.1.2)$$
for any complex number $\alpha,$ where the principal value of
$q^\alpha$ is
(usually) taken when $q \ne 0.$
In particular,  if $0<|q|<1$ and  $n=0,1,2,\ldots,$ then
$$(a;q)_{-n} = {(a;q)_\infty\over (a
q^{-n};q)_\infty}={1\o (aq^{-n};q)_n}= {(-q/a)^n\o (q/a;q)_n}\,
q^{
n\choose 2}\eqno (2.1.3)$$
and
$${1 \o (q^n;q)_k} =(q^{n+k};q)_{-k} = 0,
\qquad n\ge 1, \ k = -n,-n-1, \ldots \ . \eqno (2.1.4)$$
\par
Let $0<|q|<1, 0<|z|<1,$ and $\ n=1,2,\ldots \ .$
 Then the $q$-binomial theorem (1.2.8) gives
$${(az;q)_\infty\over (z;q)_\infty}=
\som{k=-\infty}{\infty} {(a;q)_k\over(q;q)_k}\, z^k 
=\som{j=-\infty}{\infty} {(a;q)_{j+n-1}\over(q;q)_{j+n-1}}\,
z^{j+n-1} $$
$$=\som{j=-\infty}{\infty} {(a;q)_{n-1}(aq^{n-1};q)_{j}\over
(q;q)_{n-1}(q^n;q)_{j}}\, z^{j+n-1} 
={(a;q)_{n-1}\over(q;q)_{n-1}}\, z^{n-1} 
\som{j=-\infty}{\infty} {(aq^{n-1};q)_{j}\over
(q^n;q)_{j}}\, z^{j} 
\eqno(2.1.5)$$
by a shift in the index of summation.
After replacing $a$ by $aq^{1-n} $ and then setting $q^n = b$ it
follows from the left and right sides of (2.1.5) that
$$\som{j=-\infty}{\infty} {(a;q)_{j}\over
(b;q)_{j}}\, z^{j} ={(q,b/a,az,q/az;q)_\ty \o
(b,q/a,z,b/az;q)_\ty} \eqno(2.1.6)$$
for $b=q^n$ when $0<|q|<1, 0<|z|<1,$ and $\ n=1,2,\ldots,$
where we used the compact notation
$$(a_1, a_2, \ldots, a_k;q)_\ty= (a_1;q)_\ty (a_2;q)_\ty\cdots
(a_k;q)_\ty.\eqno(2.1.7)$$
\par
By applying (2.1.3) to the terms of the series in (2.1.6) with
negative $j,$
we get 
$$\som{j=-\infty}{\infty} {(a;q)_{j}\over(b;q)_{j}}\, z^{j}
=\som{j=0}{\infty} {(a;q)_{j}\over(b;q)_{j}}\, z^{j}
+\som{j=1}{\infty} 
{(q/b;q)_{j}\over(q/a;q)_{j}}\bigl({b \o
az}\bigr)^{j},\eqno(2.1.8)$$
from which it is clear that the bilateral series on the left side
of (2.1.8)
converges in the annulus $|b/a|<|z|<1$ when $|q|<1$ and $a \ne
0.$
Since both sides of (2.1.8) are analytic functions of $b$ when
$|b|< \min (1, |az|)$ and $|z|<1,$ and since (2.1.6) holds for
 $b=q^n$ when $ 0<|q|<1,\ n=1,2,\ldots,$ it follows by analytic
continuation
that we have derived the summation formula
$$\som{j=-\infty}{\infty} {(a;q)_{j}\over
(b;q)_{j}}\, z^{j} ={(q,b/a,az,q/az;q)_\ty \o
(b,q/a,z,b/az;q)_\ty}, \qquad |b/a|<|z|<1, \ |q|<1,\eqno(2.1.9)$$
which is an extension of the $q$-binomial theorem. This formula
was stated
without proof in Ramanujan's second notebook (see Hardy [1940]
and
Berndt [1993]). The
first published proofs were given by
Hahn [1949] and M. Jackson [1950].  The above proof is
essentially the reverse
order of Ismail's [1977] proof, which first reduces the proof of
(2.1.9)
to the case when $b=q^n,$ where $n$ is a positive integer,
and then proves this case by
using a shift in the index of summation to obtain a series that
is summable by the $q$-binomial theorem.
 In the next section,
when we introduce the $_r\psi_s$ notation for general bilateral
basic
hypergeometric series, it will be seen that  (2.1.9) gives the
sum of a $_1\psi_1$ series.  Hence, (2.1.9)  is sometimes called
Ramanujan's
 $_1\psi_1$ summation formula.
\par
In BHS  a proof due to Andrews and Askey [1978] is given, which
first
considers the sum of the series in (2.1.9) as a function of $b$,
 say $f(b),$ and then shows that this function satisfies the
functional 
equation
$$f(b)={1-b/a\o (1-b)(1-b/az)}\, f(bq).
\eqno(2.1.10)$$
Iterating (2.1.10) $n-1$ times gives
$$f(b)={(b/a;q)_n\o (b,b/az;q)_n}\, f(bq^n).
\eqno (2.1.11)$$
Since $f(b)$ is analytic for $|b|<\min (1,|az|),$ we may let
$n\to\infty$
in (2.1.11) to obtain
$$f(b)={(b/a;q)_\infty\o (b,b/az;q)_\infty}\,f(0).
\eqno (2.1.12)$$
To calculate $f(0),$ it suffices to set $b=q$ in (2.1.12) and
observe that
$$f(0)={(q,q/az;q)_\infty\o (q/a;q)_\infty}\, f(q) =
{(q,q/az,az;q)_\infty\o (q/a,z;q)_\infty},\eqno (2.1.13)$$
since
$$f(q)=\sum^\infty_{n=0} {(a;q)_n\o (q;q)_n}\,z^n =
{(az;q)_\infty\o (z;q)_\infty} $$
by the $q$-binomial theorem.  Combining (2.1.12) and 
(2.1.13) gives (2.1.9).
\par
Just as in the third derivation given in \S 1.2 for the
$q$-binomial theorem, 
one can also derive (2.1.9) by using a functional equation and a
recurrence
relation to
calculate the coefficients $c_n$ in the Laurent
expansion
$${(az,q/az;q)_\ty \o (z,b/az;q)_\ty}
=\sum_{n=-\infty}^{\infty} c_n z^n. \eqno (2.1.14)$$
See, e.g., Venkatachaliengar [1988], whose proof is
reproduced in Berndt [1993] (also see Exercise 2.4).
\bigskip
{\bf  2.2 Jacobi's triple
product identity and theta functions.}  If we set $b=0$ in
(2.1.9),
replace $q$ and $z$ by $q^2$ and $-qz/a,$ respectively, and then
let $a \to \ty,$ 
we obtain Jacobi's [1829] {\it triple product identity }
$$\sum^\ty_{k=-\ty} q^{k^2}z^k = \l(q^2, -qz, -
q/z;q^2\r)_\ty.\eqno (2.2.1)$$
In BHS this formula is derived by using Heine's [1847] 
summation formula for a $_2\phi_1(c/ab)$ series,
which we will consider in \S 3.

Jacobi's triple product identity has many important applications.
In particular, it can be used to
express the {\it theta functions} (Whittaker and Watson [1965,
Chapter 21])
$$\eqalignno{
&\vartheta_1(x) = 2 \sum^\ty_{n=0} (-1)^n q^{(n+ 1/2)^2} \sin
(2n+1)x,&(2.2.2)\cr
&\vartheta_2(x) = 2 \sum^\ty_{n=0} q^{(n+ 1/2)^2}\cos
(2n+1)x,&(2.2.3)\cr
&\vartheta_3 (x) = 1 + 2 \sum^\ty_{n=1} q^{n^2} \cos 2
nx,&(2.2.4)\cr
&\vartheta_4 (x) = 1 + 2 \sum^\ty_{n=1} (-1)^n q^{n^2} \cos
2nx&(2.2.5)\cr}$$
in terms of infinite products. To derive the infinite product
representations,
 replace $q$ by $q^2$ in
(2.2.1) and then set $z$ equal to $qe^{2ix}, -qe^{2ix}, -e^{2ix},
e^{2ix}$, respectively, to obtain
$$\eqalignno{
&\vartheta_1(x) = 2q^{1/4} \sin x \prod^\ty_{n=1} (1-q^{2n}) (1-
2q^{2n}\cos 2x + q^{4n}),&(2.2.6)\cr
&\vartheta_2 (x) = 2q^{1/4} \cos x \prod^\ty_{n=1} (1-q^{2n})
(1+2q^{2n}\cos 2x + q^{4n}),&(2.2.7)\cr
&\vartheta_3(x) = \prod^\ty_{n=1} (1-q^{2n}) (1+2q^{2n-1} \cos 2x
+q^{4n-2}),&(2.2.8)\cr
\noalign{\hbox{and}}
&\vartheta_4(x) = \prod^\ty_{n=1} (1-q^{2n}) (1-2q^{2n-1}\cos 2x
+q^{4n-2}).&(2.2.9)\cr}$$

For other applications of (2.2.1), including a proof of the
{\it quintuple product identity} 
$$
 \sum^\infty_{n=-\infty} (-1)^n q^{n(3n-1)/2}z^{3n}(1+zq^n) 
= (q,-z,-q/z;q)_\infty (qz^2,q/z^2;q^2)_\infty, \qquad z\not=
0, \eqno(2.2.10)$$
see Berndt [1993] and the Notes to Exercise 5.6 in BHS.
\bigskip
\centerline{\bf Exercises 2}
\bigskip
\item{2.1}  Verify the identities (2.1.3) and (2.1.4), and show
that
$$\leqalignno{(aq^k;q)_{n-k} &= {(a;q)_n\over (a;q)_k},&(i)\cr
(aq^{2k};q)_{n-k} &= {(a;q)_n(aq^n;q)_k\over
(a;q)_{2k}},&(ii)\cr
(aq^{-n};q)_k &= {(a;q)_k (qa^{-1};q)_n\over
(a^{-1}q^{1-k};q)_n} \, q^{-nk},&(iii)\cr
(q^{-n};q)_k &= {(q;q)_n\over (q;q)
_{n-k}}(-1)^kq^{{\scriptstyle k\mathstrut\choose\scriptstyle
2}-nk}, &(iv)\cr
(a;q)_{n-k}& = {(a;q)_n\over (a^{-1}q^{1-n};q)_k}\;
(-qa^{-1})^kq^{{\scriptstyle k\mathstrut\choose\scriptstyle
2}-nk},&(v)\cr}$$
\noindent
when $k$ and $n$ are integers and both sides of the identity are
well-defined.
\medskip
\item{2.2}
Use series manipulations to verify that the function $f(b)$
defined in \S 2.1 
satisfies the functional equation (2.1.10).
\medskip
\item{2.3} Apply the inversion identity (1.3.10) to the
$q$-shifted factorials
in the left side of (2.1.9) to derive a bilateral extension of
the
$|q| > 1$ case of the $q$-binomial given in (1.3.11).
\medskip
\item{2.4}
Show that the function 
$$h(z) ={(az,q/az;q)_\ty \o (z,b/az;q)_\ty}$$
satisfies the functional equation $(b-aqz) \ h(qz)= q(1-z) \
h(z)$ 
when $ |b/aq|<|z|<1.$
Use this equation to derive a recurrence relation for
 the coefficients $c_n$ in (2.1.14). As in the derivation
of (1.2.17), use the
recurrence relation to calculate these coefficients, and hence
derive
Ramanujan's $_1\psi_1$ summation formula (2.1.9).
See  Venkatachaliengar [1988] and  Berndt [1993]. 
\medskip
\item{2.5} Investigate what happens if, as in Ex. 2.4,
you try to generalize Ramanujan's $_1\psi_1$ summation formula
by using a functional equation and a recurrence relation
to calculate the coefficients in the Laurent expansion 
$${(az,b/z;q)_\ty \o (z,c/z;q)_\ty}
=\sum_{n=-\infty}^{\infty} d_n \, z^n,  \qquad |c|<|z|<1.$$
\medskip
\item{2.6}
 Use functional equations in $a$ and $b$ to prove that
$$\int^\ty_0 {(-tq^b, -q^{a+1}/t;q)_\ty\o (-t, -q/t;q)_\ty} {dt\o
t} = -\log  q \; {(q,q^{a+b};q)_\infty\over(q^a,q^b;q)_\infty}$$
 when $ 0 < q < 1,\ {\rm Re}\ a > 0$ and Re $b > 0.$
See Askey and Roy [1986], and Gasper [1987].
\medskip
\item{2.7} Extend the integral formula in Ex. 2.6 by evaluating
the 
more general integral
$$\int^\ty_0 t^{c-1} {(-tq^b, -q^{a+1}/t;q)_\ty\o (-t, -
q/t;q)_\ty}\ dt $$
 when $ 0 < q < 1,\ {\rm Re}\ (a+c) > 0$ and Re $(b-c) >0.$
See Ramanujan [1915], Askey and Roy [1986], and Gasper [1987].
\bigskip
\centerline{\bf  3 Basic hypergeometric series, $q$-gamma and
$q$-beta
functions, $q$-integrals,}
\centerline{\bf   and some important summation and
transformation formulas}
\bigskip
{\bf  3.1  Basic hypergeometric series.}  So far, because of
the relatively simple $q$-series considered, we have  not needed
to
introduce a compact notation for $q$-series containing several
parameters.
Recall that the  {\it Gauss} [1813] 
{\it hypergeometric series} is (formally) defined by
$$F (a,\;b;\;c;\;z)\equiv {}_2F_1 (a,\;b;\;c;\;z)\equiv 
{}_2F_1\left[\matrix{a, b\cr
c\cr};z\right] =
\som{n=0}{\infty} {(a)_n(b)_n\over n!(c)_n 
}\;z^n, \eqno (3.1.1)$$
 where it is assumed that $c\not=0,-1,-2,...$, so that no zero
factors appear in the denominators of the terms of the series.
Gauss' series converges absolutely for $|z|<1$, and
for $|z| = 1$ when  Re\nobreak $(c-a-b) >0.$ 
 Heine [1846, 1847, 1878] introduced the series
$$\phi (\alpha, \beta, \gamma,q,z)={}_2\phi_1 (q^\alpha,
 q^\beta; q^\gamma; q,z)
\eqno (3.1.2)$$
with
$$ _2\phi_1 (a,b;c;q,z)
\equiv {}_2\phi_1\left[\matrix{a, b\cr  c\cr};q,z\right]=
\som{n=0}{\infty} {(a;q)_n(b;q)_n\over(q;q)_n(c;q)_n}\;z^n,\eqno
(3.1.3)$$
where it is assumed that $\gamma \not=-m$
and $c\not= q^{-m}$ for $m = 0, 1, \ldots\ $.
{\it Heine's series} converges absolutely for $|z|<1$ when $|q|
<1,$  
and it is a $q$-analogue of Gauss' series because, by using
(1.2.1) and
taking a formal termwise limit,
$$\lim_{q\to1^-} {}_2\phi_1 (q^\alpha, q^\beta; q^\gamma; q,z)
= {}_2F_1 (\alpha, \beta; \gamma;z).
\eqno (3.1.4)$$
In view of the base $q,$ Heine's series is
also called  the {\it basic hypergeometric series} or
$q$-{\it hypergeometric series.}  We prefer to use the $_2\phi_1
(a,b;c;q,z)$
notation instead of  Heine's $\phi (\alpha, \beta, \gamma,q,z)$
notation,
because when $0<|q| <1$ the $\alpha, \beta, $ or $\gamma \to \ty
$
 limit cases of Heine's  series
correspond to setting $a, b,$ or $c,$
 respectively, equal to zero in  (3.1.3).

The {\it (generalized)
hypergeometric
series} with $r$ numerator parameters $a_1, \ldots, a_r$ and $s$
denominator parameters $b_1, \ldots, b_s$ is defined (formally)
by
$$\eqalignno{&_rF_s (a_1,a_2,\ldots, a_r; \;b_1,\ldots
,b_s;z)\equiv
{}_rF_s\left[\matrix{a_1, a_2,\ldots, a_r\cr b_1,\ldots, b_s
\cr};z\right]\cr
&= \som{n=0}{\infty} {(a_1)_n(a_2)_n\cdots
(a_r)_n\over n!(b_1)_n \cdots (b
_s)_n}\;z^n, &(3.1.5)\cr}$$
 and an $_r\phi_s$ {\it
basic hypergeometric series} is defined by
$$\eqalign{&_r\phi_s (a_1, a_2, \ldots, a_r; b_1,\ldots,b_s;q,z)
\equiv {}_r\phi_s\left[\matrix{a_1, a_2,\ldots , a_r \cr
b_1,\ldots ,b_s\cr};q,z\right]\cr
&= \som{n=0}{\infty} {(a_1,a_2,\ldots, a_r;q)_n
\over(q,b_1,\ldots, b_s;q)_n} \left[(-1)^n
q^{n\choose 2}\right]^{1+s-r}z^n\cr}\eqno (3.1.6)$$
where ${n\choose 2}=n(n-1)/2$ and, analogous to (2.1.7),
 we employed the compact notation
 $(a_1,a_2,\ldots, a_r;q)_n = (a_1;q)_n
(a_2;q)_n\cdots(a_r;q)_n.$ 
As in (3.1.4), $_r\phi_s$ is a $q$-analogue
of $_rF_s$ and we have (formally) 
$$\lim_{q\to1^-}
{}_r\phi_s\left[\matrix{q^{a_1}, q^{a_2},\ldots , q^{a_r} \cr
{q^b_1},\ldots ,q^{b_s}\cr};q,{(q-1)^{1+s-r}z}\right]
= {}_rF_s\left[\matrix{a_1, a_2,\ldots, a_r\cr b_1,\ldots, b_s
};z\right]. \eqno (3.1.7)$$  

In terms of these notations, the binomial theorem (1.1.3)
and the $q$-binomial theorem (1.2.8) may be written in the forms
$${}_2F_1 (a,c;c;z) = {}_1F_0 (a;\hbox{---};z) = 
(1-z)^{-a},\quad |z|<1, \eqno (3.1.8)$$
 and
$${}_2\phi_1 (a,c;c; q,z) = 
{}_1\phi_0 (a;\hbox{---} ; q,z) =
{(az;q)_\infty\over (z;q)_\infty},\quad |z|<1,\ |q|
<1,\eqno
(3.1.9)$$
where a dash is used to indicate the absence of either numerator
(when $r=0$) or denominator (when $s=0$) parameters. 
Many other important special cases are
 considered in BHS, its references,
 and in what follows.

It is assumed in (3.1.5) and (3.1.6) that the parameters $b_1,
\ldots ,b_s$ are such that the denominator factors in the terms
of the series are never zero, and in (3.1.6)
it is assumed that $q\ne 0$ when $r> s+1$.  From
$$(-m)_n = \left(q^{-m};q\right)_n = 0,\quad
n=m+1, m+2, \ldots,\eqno (3.1.10)$$
\noindent 
we see that an $_rF_s$ series terminates if one of its numerator
parameters
is zero or a negative integer, and an $_r\phi_s$ series
terminates if one of its numerator parameters is of the form
$q^{-m}$ with $m=0,1,2,\ldots$, and $q\ne 0$. 
Unless stated otherwise, when working with nonterminating basic
hypergeometric series, we
will (for simplicity)  assume that $|q| <1$ and that the
parameters and
variables are such that the series converge absolutely.  As in
our 
derivation of (1.3.11),
if $|q|>1$ then we can use (1.3.10) to perform an inversion with
respect to the base to transform the series (3.1.6) into a series
in 
 base $q^{-1},$ with $|q^{-1}| <1$ (see Exercise 3.2).

The ratio test shows
that an $_rF_s$ series converges absolutely for
all $z$ if $r\leq s$, and for $|z| <1$ if $r=s+1$.  By Raabe's
test
this series also
converges absolutely for $|z| = 1$ if $r=s+1$ and
Re $[b_1+\cdots +b_s -(a_1 +\cdots + a_r)] >0.$  If $r>s+1$ and
$z\not= 0$ or $r=s+1$ and $|z| >1$, then this series
diverges, unless it terminates.
Similarly, if $0 < |q| <1$, then the $_r\phi_s$ series converges
absolutely
for all $z$ if $r\leq s$ and for $|z| <1$ if $r = s+1$.
It converges absolutely if $|q| >1$ and
$|z| < |b_1b_2\cdots b_s|/|a_1a_2\cdots a_r|$. It diverges for
$z\not= 0$ if $0 < |q|
<1$ and $r>s+1$, and if $|q| >1$ and $|z| > |b_1b_2\cdots b_s|/|
a_1 a_2 \cdots a_r|$,
unless it terminates.  The $_rF_s$ and
$_r\phi_s$ notations are also used
for the sums of these series inside the circle of convergence and
for their analytic continuations (called {\it hypergeometric
functions} and {\it basic hypergeometric functions},
respectively) outside the circle of convergence.

The series (3.1.6) has the property that if we
replace $z$ by $z/a_r$ and let $a_r\rightarrow \ty$ (called a
{\it confluence
process}), then
we obtain a series that is of the form (3.1.6) with $r$
replaced by $r-1$. 
This is not the case for the
$_r\phi_s$ series defined without the factors
$[(-1)^nq^{ n \choose 2}]^{1+s-r}$ in  Bailey
[1935] and Slater [1966].
It is because we need a notation that includes such limit cases 
as special cases that we have
chosen to use the definition (3.1.6).  Also, there is no loss
in generality because the Bailey and Slater series can be
obtained
from the $r=s+1$ case of (3.1.6) by choosing $s$ sufficiently
large and setting some of the parameters equal to zero.

Notice that if we  denote the terms of the series (3.1.5)
and (3.1.6) which contain $z^n$ by $u_n$
and $v_n$, respectively, then ${u_{n+1}\over u_n}$ is a rational
function 
of $n$ and
${v_{n+1}\over v_n}$ is a rational function of $q^n$.
Conversely, any rational function of $n$  can be written in the
form of 
${u_{n+1}\over u_n},$ and any rational function of $q^n$ can be
written in the form of 
${v_{n+1}\over v_n}.$ Hence, we have the characterization that if
$\sum^\ty_{n=0} u_n$ and $\sum^\ty_{n=0} v_n$ are
series with $u_0 = v_0 =1$ such that $u_{n+1}/u_n$ is a
rational function of $n$ and $v_{n+1}/v_n$ is a rational function
of $q^n$, then these series are of the forms (3.1.5)
and (3.1.6), respectively.  This  characterization  provides
another reason
why we defined $_r\phi_s$ in (3.1.6) with the 
$[(-1)^nq^{ n \choose 2}]^{1+s-r}$ factor.

Generalizing the bilateral $q$-series that we considered in \S
2.1, 
we define the {\it  bilateral basic hypergeometric series} 
in base $q$
with $r$
numerator and
$s$ denominator parameters  by
$$\eqalign{
_r\psi_s (z) &\equiv \;_r\psi_s \left[\matrix
{a_1, a_2,\dots,a_r \cr
b_1, b_2,\dots,b_s \cr};q,z\right] \cr
&= \sum^\infty_{n=-\infty}
{(a_1,a_2,\dots,a_r;q)_n\over (b_1,b_2,\dots,b_s;q)_n}\,(-1)
^{(s-r)n}
q^{(s-r){n\choose 2}}\, z^n,
\cr}\eqno (3.1.11)$$
where it is assumed that $q$, $z$ and the parameters are
such that each
term of the series is well-defined (i.e., the denominator factors
are never
zero, $q\ne0$ if $s<r$, and $z\ne0$ if negative powers of $z$
occur). 
A bilateral basic hypergeometric series may be characterized as
being a series
$\sum^\infty_{n=-\infty} v_n$ such that $v_0 = 1$ and
$v_{n+1}/v_n$ is a rational function of $q^n$. 

Employing the  $_1\psi_1$ notation,  Ramanujan's summation
formula (2.1.9)
 may be stated  in the form 
$$ _1\psi_1(a; b; q, z)  ={(q,b/a,az,q/az;q)_\ty \o
(b,q/a,z,b/az;q)_\ty}, \qquad |b/a|<|z|<1, \
|q|<1.\eqno(3.1.12)$$

As in our derivation of (2.1.8), to
determine when $_r\psi_s$  series converge we first apply
(1.3.10) to the terms with negative $n$ to obtain the
decomposition
$$\eqalign{
_r\psi_s(z) &= \sum^\infty_{n=0}
{(a_1,a_2,\dots,a_r;q)_n\over (b_1,b_2,\dots,b_s;q)_n}\,(-1)
^{(s-r)n}q^{(s-r){n\choose 2}}\,z^n \cr
&\quad + \sum^\infty_{n=1}
{(q/b_1,q/b_2,\dots,q/b_s;q)_n\over (q/a_1,q/a_2,\dots,q/a_r;q)_
n}\,
\left({b_1\cdots b_s\over a_1\cdots a_rz}\right)^n.\cr}
\eqno (3.1.13)$$
Set $R = |(b_1\cdots b_s)/(a_1\cdots a_r)|.$ 
From (3.1.13) it is obvious that if $r>s$ and $0<|q|<1$, 
then the first series on the right side of
(3.1.13)
diverges for $z\not= 0$; if $r>s$ and $|q| > 1$, then the first
series converges for $|z| < R$ and the second series converges
for all $z \ne 0$. When $r < s$ and $|q| < 1$ the first
series converges for all $z$, but the second series
converges only when
$|z| > R$. If $r<s$
and $|q|>1$, the second series diverges for all $z\ne0$. If
$r=s$ and $|q|<1$, the first series converges
when $|z|<1$
and the second when $|z| >R$; while
if $|q|>1,$ then the second series converges when $|z|>1$ and the
first
when $|z| < R$. In particular,
if $r=s$ and $|q|<1$, which is the most important case, 
the region of convergence of the series $ _r\psi_r(z)$
is the annulus
$\left|(b_1\cdots b_r)/( a_1\cdots a_r)\right| < |z|<1.$
\bigskip
{\bf  3.2 $q$-gamma and $q$-beta functions.}  Analogous to 
Gauss' infinite product representation for the gamma function
$$\G(z)=z^{-1} \prod^\ty_{k=1} [(1+1/k)^z (1+z/k)^{-1}],\eqno
(3.2.1)$$
the $q$-{\it gamma  function} $\G_q(z)$ is defined as in 
Thomae [1869] and Jackson [1904] by
$$\Gamma_q(z) = {(q;q)_\infty\over (q^z;q)_\infty}\,
(1-q)^{1-z}\;,\;0<q<1.\eqno (3.2.2)$$
When
$z = n+1$ is a positive integer, this definition reduces to
$$\Gamma_q(n+1) = 1 (1 +q) (1+q+q^2) \cdots (1 + q + \cdots
+q^{n-1}),\eqno (3.2.3)$$
\noindent which tends to $n!$ as $q\to 1^-$.
Thus $\Gamma_q(n+1) \rightarrow \Gamma (n+1) = n!$ as
$q\rightarrow 1^-$. One can extend the definition of $\G_q (z)$ 
 to $|q| < 1$ by using the principal values of $q^z$
and $(1-q)^{1-z}$ in (3.2.2). In BHS it is shown that
$$\lim_{q\to1^-}\Gamma_q(z) = \Gamma (z)\eqno (3.2.4)$$
and that $\Gamma_q(z)$ satisfies the 
functional equation
$$f(z+1) = {1-q^z\over 1-q} f(z),\quad f(1) = 1, \eqno(3.2.5)$$
which is a $q$-analogue of the well-known functional equation for
the gamma 
function
$$\Gamma (z+1) = z \ \Gamma (z), \quad \Gamma (1) =
1.\eqno(3.2.6)$$

Analogous to the definition of the
 the {\it beta function} 
$$B(x,y)= {\Gamma (x)\Gamma (y)\over \Gamma
(x+y)},\eqno(3.2.7)$$
\noindent  the $q$-{\it beta function} is defined by
$$B_q(
x,y) = {\Gamma_q(x)\Gamma_q(y)\over
\Gamma_q(x+y)},\eqno(3.2.8)$$
\noindent which tends to $B(x,y)$ as $q\rightarrow 1^-$.  From
(3.2.2) and (1.2.8),
$$\eqalignno{B_q(x,y) &= (1-q)
{(q,q^{x+y};q)_\infty\over(q^x,q^y;q)_\infty}\cr
&= (1-q) {(q;q)_\infty\over
(q^y;q)_\infty} \som{n=0}{\infty}{(q^y;q)_n\over
(q;q)_n}\, q^{nx}\cr
&= (1-q) \som{n=0}{\infty}
{(q^{n+1};q)_\infty\over(q^{n+y};q)_\infty}\, q^{nx},\quad
{\rm Re}\;x,\; {\rm Re}\;y>0.&(3.2.9)\cr}$$
\noindent This expansion will be utilized in \S 3.3 
to derive a $q$-integral representation for
$B_q(x,y)$.
\bigskip
{\bf  3.3   $q$-integrals.} Thomae [1869, 1870] defined
the $q$-{\it integral} on the interval $[0, 1]$  by
$$\int_{0}^{1} f(t)\ d_qt = (1-q)
\som{n=0}{\infty} f(q^n)q^n.\eqno(3.3.1)$$
The right side of (3.3.1) corresponds to using a Riemann sum with
partition points $t_n = q^n, n=0,1,2,\ldots \ .$\ \
\noindent Jackson [1910b] extended this to the interval $[a, b]$
via 
$$\int_{a}^{b} f(t)\ d_qt = \int_{0}^{b}f(t)\ d_qt -
\int_{0}^{a}f(t)\ d_qt,\eqno (3.3.2)$$
where
$$\int_{0}^{a}f(t)\ d_qt = a(1-q) \som{n=0}{\infty}f(a
q^n)\, q^n. \eqno (3.3.3)$$
He also defined an integral on ($0, \infty$) by
$$\int_{0}^{\infty} f(t)\ d_qt = (1-q)
\som{n=-\infty}{\infty}f(q^n)q^n.\eqno (3.3.4)$$
On the interval ($-\ty, \infty$) the 
{\it bilateral} $q$-{\it integral} is defined by
$$\int_{-\infty}^{\infty} f(t)\ d_qt = (1-q) \som
{n=-\infty}{\infty} \left[f(q^n) + f(-q^n)\right]q^n.\eqno
(3.3.5)$$

When $f$ is continuous on $[0,a]$, it can be shown that
$$\lim_{q\to1} \int_{0}^{a}f(t)\ d_qt = \int_{0}^{a}f(t)\ dt
\eqno(3.3.6)$$
\noindent and that a similar limit holds for (3.3.4)
and (3.3.5) for suitable functions.  
By employing the $q$-integral definition (3.3.1), the 
series expansion (3.2.9) for the $q$-beta function can be
rewritten in the form
$$B_q(x,y) = \int_{0}^{1} t^{x-1}
{(tq;q)_\infty\over(tq^y;q)_\infty}\ d_qt, \quad {\rm Re}\;x
>0,\quad y \ne 0, -1, -2, \ldots, \eqno(3.3.7)$$
\noindent which, as $q\rightarrow 1^-$, tends to 
the beta function integral
$$B(x,y) = \int_{0}^{1} t^{x-1}(1-t)^{y-1}\ dt,\qquad
{\rm Re}\; x,\; {\rm Re}\; y >0.\eqno (3.3.8)$$
\noindent 

We can derive a $q$-integral representation for $_2\phi_1
(q^a,q^b;q^c;q,z)$
by using (3.3.7) to get
$${(q^b;q)_n\over(q^c;q)_n}=
{\Gamma_q(c)\over\Gamma_q(b)\Gamma_q(c-b)} 
\int_{0}^{1} t^{b+n-1}
{(tq;q)_\infty\over(tq^{c-b};q)_\infty}\ d_qt \eqno (3.3.9)$$
for $n \ge 0, \; {\rm Re}\;b>0, \;c-b \ne 0, -1, -2, \ldots,$ 
which leads to
$$_2\phi_1 (q^a,q^b;q^c;q,z)
= {\Gamma_q(c)\over\Gamma_q(b)\Gamma_q(c-b)} 
\som{n=0}{\infty}{(q^a;q)_n\over(q;q)_n} z^n
\int_{0}^{1} t^{b+n-1}
{(tq;q)_\infty\over(tq^{c-b};q)_\infty}\ d_qt$$
$$=
{\Gamma_q(c)\over\Gamma_q(b)\Gamma_q(c-b)} 
\int_{0}^{1} t^{b-1}
{(tq;q)_\infty\over(tq^{c-b};q)_\infty} 
{}_1\phi_0 (q^a;\hbox{---} ; q,tz) \ d_qt \eqno(3.3.10)$$
when $|z|<1,$ ${\rm Re}\;b>0,$ and $\; c-b \ne 0, -1, -2, \ldots
\ .$
Employing (3.1.9) to sum the ${}_1\phi_0$ series in 
(3.3.10) gives Thomae's [1869] $q$-integral
$$_2\phi_1 (q^a,q^b;q^c;q,z)
= {\Gamma_q(c)\over\Gamma_q(b)\Gamma_q(c-b)} \int_{0}^{1} t^{b-1}
{(tzq^a,tq;q)_\infty\over (tz, tq^{c-b};q)_\infty}\
d_qt,\eqno(3.3.11)$$
where $|z|<1, {\rm Re}\; b>0,$ 
and $c-b \ne 0, -1, -2, \ldots \ $.
This is a $q$-analogue of Euler's integral
representation
$$_2F_1 (a,b;c;z) = {\Gamma (c)\over \Gamma (b)
 \Gamma (c-b)}\int_{0}^{1} t^{b-1}(1-t)^{c-b-1}
(1-tz)^{-a}\ dt,\eqno (3.3.12)$$
\noindent where $|\arg (1-z)| <\pi$ and Re $c > \; {\rm Re}\; b
>0$.
\bigskip
{\bf 3.4 Heine's $_2\phi_1$ transformation and summation
formulas.}
By employing Thomae's definition (3.3.1)  of a  $q$-integral  to
rewrite
the integral in (3.3.11) as a series, we obtain 
 Heine's [1847], [1878] $_2\phi_1$ transformation
formula
$$_2\phi_1 (q^a,q^b;q^c;q,z) = {(q^b,q^az;q)_\infty\over
(q^c,z;q)_\infty}\;_2\phi_1 (q^{c-b}, z; q^az;q,q^b),\eqno
(3.4.1)$$
or, equivalently, on replacing $q^a, q^b,q^c$ by $a, b, c,$
respectively,
$$_2\phi_1 (a,b;c;q,z) = {(b,az;q)_\infty\over
(c,z;q)_\infty}\;_2\phi_1 (c/b, z; az;q,b),\eqno (3.4.2)$$
\noindent provided $|z| <1$ and $|b| <1$.   
Thomae [1869] had derived his integral representation (3.3.11) by
rewriting
Heine's formula (3.4.2) in the form (3.3.11), the reverse of our
derivation
of (3.4.2).  Formula  (3.4.2) is derived in BHS by
using the $q$-binomial theorem and series manipulations of double
series.

One particularly attractive feature of  (3.4.2) is that it can be
iterated
to produce the following chain of transformation formulas
$$\eqalignno{_2\phi_1 (a,b;c;q,z) &=
{(b,az;q)_\infty\over(c,z;q)_
\infty} \;_2\phi_1(c/b,z;az;q,b)&(3.4.3)\cr
&= {(c/b, bz;q)_\infty\over
(c,z;q)_\infty}\;_2\phi_1(abz/c,b;bz;q,c/b)&(3.4.4)\cr
&= {(abz/c;q)_\infty\over (z;q)_\infty}\;_2\phi_1
(c/a,c/b;c;q,abz/c),&(3.4.5)\cr}$$
provided the series converge.
In particular, (3.4.5) and the left side of (3.4.3) yield Heine's
transformation formula
$$_2\phi_1 (a,b;c;q,z) = {(abz/c;q)_\infty\over
(z;q)_\infty}\;_2\phi_1 (c/a, c/b; c; q,abz/c)\eqno (3.4.6)$$
when $|z|<1$ and $|abz/c|<1,$
which is a $q$-analogue of Euler's transformation formula for
$_2F_1$ series
$$_2F_1 (a,b;c;z) = (1-z)^{c-a-b} \;_2F_1 (c-a, c-b; c;z), \qquad
|z|<1.\eqno
(3.4.7)$$
Formulas (3.4.2) -- (3.4.6) and  (3.3.11)
can be used to analytically continue  $_2\phi_1 (a,b;c;q,z)$ 
functions to regions outside the unit disk $|z|<1.$ 
\noindent

Another important application of (3.4.2) is that in the case
 $z=c/ab$ it reduces to
$$_2\phi_1 (a,b;c;q,c/ab) =
{(b,c/b;q)_\infty\over(c,c/ab;q)_\infty}\; _1\phi_0 (c/ab;
\hbox{---};
q,b), \qquad |c/ab|<1,\  |b|<1, \eqno (3.4.8)$$
where the $_1\phi_0$ can be summed via the $q$-binomial
theorem to yield Heine's [1847] summation formula for a
$_2\phi_1$ series
$$_2\phi_1 (a,b;c;q,c/ab) = {(c/a, c/b;q)_\infty\over(c,
c/ab;q)_\infty}, \eqno (3.4.9)$$
which, by analytic continuation, holds when $|c/ab|<1.$ 
This formula is Heine's [1847] $q$-analogue
 of Gauss' [1813] famous summation formula
$$F(a,b;c;1) = {\Gamma (c)\Gamma (c-a-b)\over \Gamma
(c-a)\Gamma(c-b)},\quad {\rm Re}\; (c-a-b) >0.\eqno (3.4.10)$$
In \S 1.6 of BHS it is shown that the Jacobi triple product
identity (2.2.1)
can be easily derived from Heine's summation formula (3.4.9).

For the terminating case when $b=q^{-n}$, (3.4.9) reduces to
$$_2\phi_1 ( a,q^{-n};c;q,cq^n/a) = {(c/a
;q)_n\over (c;q)_n},\eqno(3.4.11)$$
 which, by using the inversion identity (1.3.10) or
by changing the order of summation, is equivalent to
$$_2\phi_1 (a,q^{-n};c;q,q) = {(c/a;q)_n\over (c;q)_n}\, 
a^n.\eqno(3.4.12)$$
Formulas (3.4.11) and (3.4.12) are  equivalent  to (1.3.3),
and they are $q$-analogues of the Chu--Vandermonde 
  summation formula
$$F(a,-n;c;1) = {(c-a)_n\over (c)_n},\quad n = 0,1, \ldots \
.\eqno
(3.4.13)$$

In \S 1.5 of BHS it is shown that (3.4.11) and (3.4.12)  can be 
used to derive  Sears' [1951b] transformation formula
$$\eqalignno{&_2\phi_1 (q^{-n}, b;c;q
,z)\cr&= {(c/b;q)_n\over (c;q)_n} \left({bz\over q}\right)^n
\;_3\phi_2(q^{-n}, q/z, c^{-1}q^{1-n};bc^{-1}q^{1-n},
0;q,q)& (3.4.14)\cr}$$
and Jackson's [1910a]
transformation formula
$$\eqalignno{_2\phi_1 (a,b;c;q,z) &=
{(az;q)_\infty\over(z;q)_\infty} \som{k=0}{\infty}
{(a,c/b;q)_k\over (q,c,a
z;q)_k}(-bz)^kq^{\scriptstyle k\mathstrut \choose\scriptstyle
2}\cr&= {(az;q)_\infty\over
(z;q)_\infty} \; _2\phi_2 (a, c/b;c,az;q,bz),&(3.4.15)\cr}$$
which is a $q$-analogue of the Pfaff-Kummer
transformation formula
$$_2F_1 (a,b;c;z) = (1-z)^{-a} \; _2F_1 (a,c-b
; c;z/(z-1)).\eqno(3.4.16)$$
\bigskip
{\bf 3.5  $q$-analogue of the Pfaff--Saalsch\"utz summation
formula.} Observe
that, since 
$${(abz/c;q)_\infty\over (z;q)_\infty} =
\som{k=0}{\infty}{(ab/c;q)_k\over (q;q)_k}\,  z^k$$
by the $q$-binomial theorem, the right side of (3.4.6) equals
$$\som{k=0}{\infty} \som{m=0}{\infty} {(ab/c;q)_k
(c/a,c/b;q)_m\over (q;q)_k (q,c;q)_m} \left({ab
\over c}\right)^mz^{k+m} ,$$
\noindent and hence, by equating the coefficients of $z^n$ on
both
sides of (3.4.6), 
$$\som{j=0}{n} {(q^{-n}, c/a, c/b;q)_j\over
(q,c,cq^{1-n}/ab;q)_j} \, q^j = {(a,b;q)_n\over (c,ab/c;q)_n}.$$
After replacing $a, b$ by $c/a,
 c/b$, respectively, this yields Jackson's [1910a] summation
formula
for a  terminating $_3\phi_2$ series
$$_3\phi_2 (a,b,q^{-n}; c, abc^{-1}q^{1-n};q,q) =
{(c/a, c/b;q)_n\over (c, c/ab;q)_n},\quad n =  0,1,
\ldots\ . \eqno (3.5.1)$$
By
replacing $a,b,c$ in (3.5.1) by $q^a, q^b,q^c$, respectively, and
letting $q\rightarrow 1$ one obtains the 
 sum of a  terminating $_3\phi_2$ series
$$_3F_2 (a,b,-n;c, 1+a+b-c-n;1) = {(c-a)_n(c-b)_n\over (c)_n
(c-a-b)_n},\quad n = 0, 1, \ldots,\eqno (3.5.2)$$
which was discovered by Pfaff [1797] and rediscovered by 
Saalsch\"utz [1890]. Since (3.5.2) is usually
called the {\it Saalsch\"utz formula} or the {\it 
Pfaff-Saalsch\"utz formula}, Jackson's formula
(3.5.1) is usually called the $q$-{\it Saalsch\"utz }
or the $q$-{\it Pfaff-Saalsch\"utz formula}.
It should be noted that letting $a\rightarrow \infty$
in (3.5.1) gives (3.4.11), while letting $a\rightarrow 0$
 gives (3.4.12).

The  Pfaff-Saalsch\"utz $_3F_2 (1)$ series is said to be a {\it
balanced}
series (or {\it
Saalsch\"utzian})
because $z = 1$ and the sum of its denominator parameters equals
the
sum of its numerator parameters plus 1.  Analogously, the 
$q$-Pfaff-Saalsch\"utz $_3\phi_2 (q) $
 series is also called  {\it balanced}
because $z = q$ and 
the product of its denominator parameters equals the
product of its numerator parameters times $q.$  
In general, as defined in BHS,
an $_{r+1}F_r$ series is called $k$-{\it balanced} if $b_1+b_2
+\cdots +b_r = k + a_1 + a_2 + \cdots +a_{r+1}$ and $z = 1$; a
1-balanced series is called {\it balanced} (or {\it
Saalsch\"utzian}). Analogously, an $_{r+1
}\phi_r$ series is called $k$-{\it balanced} if $b_1b_2\cdots b_r
=
q^k a_1a_2\cdots a_{r+1}$ and $z = q$,
and a 1-balanced series is called {\it balanced} (or {\it
Saalsch\"utzian}). In general
an $_{r+1}F_r(z)$ series is called {\it balanced} if $b_1+b_2
+\cdots +b_r = 1 + a_1 + a_2 + \cdots +a_{r+1}$ and $z = 1$, 
and an $_{r+1
}\phi_r(z)$ series is called {\it balanced} if $b_1b_2\cdots b_r
=
q a_1a_2\cdots a_{r+1}$ and $z = q$.  We will
 encounter several balanced series in \S 4.
 \bigskip
\centerline{\bf  Exercises 3}
\bigskip
\item{3.1} Verify the convergence conditions stated for
$_rF_s, \; {}_r\phi_s, $ and ${}_r\psi_s$ series in \S 3.1.
\medskip
\item{3.2} Use (1.3.10) to derive the (formal) inversion formula
$${}_r\phi_s  \l[\matrix{a_1, \ldots, a_r\cr
\noalign{\smallskip}
b_1, \ldots, b_s\cr} ;q,z\r] = \sum^\ty_{n=0} {(a^{-1}_1, 
\ldots, a^{-1}_r; q^{-1})_n\o (q^{-
1}, b^{-1}_1, \ldots, b_s^{-1};q^{-1})_n} \l({a_1 \cdots
a_rz\over b_1 \cdots b_s q}\r)^n.$$
\item{3.3} Show that if $z$ is
replaced  by $z/a_r$ in the series (3.1.6), then on letting 
$a_r\rightarrow \ty$ 
one obtains a series that is of the form (3.1.6) with $r$
replaced by $r-1$. 
Show that this is not the case if the
$_r\phi_s$ series in  (3.1.6) is defined without the factors
$[(-1)^nq^{ n \choose 2}]^{1+s-r}.$ 
\medskip
\item{3.4}  Prove the limit relation (3.2.4) and the functional
equation
(3.2.5).
\medskip
\item{3.5} Derive the transformation formulas (3.4.14) and
(3.4.15).
Show that (3.4.15) is a $q$-analogue of (3.4.16).

\medskip
\item{3.6} Derive the generating function
$$\som{n=0}{\infty} {H_n(x|q)\over (q;q)_n} \, t^n =
{1\over(te^{i\theta}, te^{-i\theta};q)_\infty},\quad  |t|
<1,$$
for the {\it continuous} $q$-{\it Hermite
polynomials} defined  by
$$H_n(x|q) = \som{k=0}{n} {(q;q)_n\over (q;q)_k
(q;q)_{n-k}}\, e^{i(n-2k)\theta},$$
\noindent where $x =
\cos\theta.$ 
See Rogers [1894] and Askey and Ismail [1983].
\medskip\item{3.7}   Derive the generating function
$$\som{n=0}{\infty} C_n (x;\beta |q)\, t^n = {(\beta
te^{i\theta}, \beta te^{-i\theta};q)_\infty\over
(te^{i\theta},te^{-i\theta};q)_\infty},\; \ |t| <1,$$
for the {\it continuous} $q$-{\it ultraspherical
polynomials} defined by
$$C_n(x;\beta |q) = \som{k=0}{n} {(\beta ;q)_k (\beta
;q)_{n-k}\over (q;q)_k (q;q)_{n-k}}\,  e^{i(n-2k)\theta},$$
\noindent where $x = \cos\theta$. 
See Rogers [1895] and Askey and Ismail [1983].
\medskip
\item{3.8} The {\it little} $q$-{\it Jacobi polynomials} are
defined in Andrews and Askey [1977] by
$$p_n (x; a,b;q) = \;_2\phi_1 (q^{-n},abq^{n+1}; aq;q,qx).$$
\noindent Prove that they satisfy the orthogonality
relation
$$\eqalign{
&\som{j=0}{\infty} {(bq;q)_j\over (q;q)_j} (aq)^j\,
p_n(q^j;a,b;q)\, 
p_m(q^j;a,b;q)\cr&= \cases{0,& if $m \not= n$,\cr
{\ds (q,bq;q)_n
(1-abq)(aq)^n\over\ds (aq,abq;q)_n (1-abq^{2n+1})}  {\ds
(abq^2;q)_\infty\over\ds (aq;q)_\infty},& if $m=n$.\cr}\cr}$$
\medskip
\item{3.9} Show that the little $q$-Jacobi polynomials defined in
Ex. 3.8
satisfy the {\it connection coefficient formula}
$$p_n(x;c,d;q) = \som{k=0}{n} a_{k,n}\ p_k(x;a,b;q)$$
\noindent  with
$$a_{k,n} = {(q^{-n},aq,cdq^{n+1};q)_k\over
(q,cq,abq^{k+1};q)_k}\;_3\phi_2\left[\matrix{q^{k-n}, 
cdq^{n+k+1}, aq^{k+1}\cr
cq^{k+1}, abq^{2k+2}\cr};q,q\right].$$
\bigskip
\centerline{\bf  4 Summation, transformation, and expansion
formulas,}
\centerline{\bf  integral representations, and  applications}
\bigskip
{\bf  4.1  Finite differences, bibasic and very-well-poised
series.}  Let
 $n, m  =  0, \pm 1, \pm 2, \ldots .$ 
Another method
to derive summation formulas is to use the fact that if a finite
difference
operator  $\D$ is defined for any sequence $\{u_k\}$ (of real or
complex
numbers), by
$\D u_k = u_k - u_{k-1},$ then
$$\sum^n_{k=m} \D u_k = u_n - u_{m-1},\eqno (4.1.1)$$
where  we employed the convention of defining
$$\eqalignno{
&\sum^n_{k=m} a_k = \cases{a_m + a_{m+1} + \cdots + a_n,& $m\le
n$,\cr
&\cr
0,&$m=n+1$,\cr
&\cr
-(a_{n+1} + a_{n+2} + \cdots + a_{m-1}),& $m\ge n + 2$,\cr}\cr
&&(4.1.2)\cr}$$
for any sequence $\{a_k\}.$  
An excellent way to motivate this definition is
to (formally) let
$$u_n =  \sum^n_{k=-\ty} a_k$$
and then observe that, by 
cancelling the $a_k$'s that appear in both series in the
difference
$$u_n  - u_{m-1} = \sum^n_{k=-\ty}a_k - \sum^m_{k=-\ty} a_k ,$$
we  obtain (4.1.1) with $\sum^n_{k=m} a_k$ defined by (4.1.2).

If, as in Gasper [1989a], we let
$$s_k = {(ap, bp;p)_k (cq, aq/bc;q)_k\o (q, aq/b;q)_k (ap/c,
bcp;p)_k}, \eqno (4.1.3)$$
where $q$ and $p$ are independent bases, 
we obtain the factorization
$$\eqalignno{
\D s_k &= s_k - s_{k-1}\cr
&= {(ap, bp;p)_{k-1} (cq,aq/bc;q)_{k-1}\o (q, aq/b;q)_k
(ap/c, bcp;p)_k}\cr
&\quad\cdot \l\{ (1-ap^k) (1-bp^k) (1-cq^k) (1-aq^k/bc)\r.\cr
&\qquad \l.- (1-q^k)(1-aq^k/b) (1-ap^k/c) (1-bcp^k)\r\}\cr
&= {(1-ap^kq^k) (1-bp^k/q^k)\o (1-a) (1-b)
}{(a,b;p)_k (c,a/bc;q)_k\o (q, aq/b;q)_k
(ap/c,
bcp;p)_k}\,q^k,&(4.1.4)\cr}$$
which is equivalent to the easily checked factorization
$$\eqalignno{
&(1-a)(1-b)(1-c)(1-ad^2/bc) - (1-d)(1-ad/b)(1-ad/c)(1-bc/d)\cr
&=(1-ad)(1-b/d)(1-ad/bc)(d-c).&(4.1.5)\cr}$$
Since $ s_k =0$ when $k \le -1,$ (4.1.4) gives the 
{\it indefinite bibasic summation formula}
$$\eqalignno{
&\sum^n_{k=0} {(1-ap^kq^k) (1-bp^kq^{-k})\o (1-a)(1-b)}
{(a,b;p)_k(c,a/bc;q)_k\o (q, aq/b;q)_k (ap/c, bcp;p)_k} q^k\cr
&= {(ap,bp;p)_n(cq,aq/bc;q)_n\o (q,aq/b;q)_n (ap/c,
bcp;p)_n}&(4.1.6)\cr}$$
for $n = 0, 1, \ldots\ $.  
In particular, since $\left(q^{1-n};q\right)_n=0$ unless $n=0$,
when $p=q$ and $c=q^{-n}, \ n = 0, 1, \ldots\,, $ 
(4.1.6) reduces to the $_6\phi_5$ summation formula

$$_6\phi_5\left[\matrix{a, qa^\half, -qa^\half, b, aq^n/b,
q^{-n}\cr
\noalign{\smallskip}
a^\half, -a^\half, aq/b, bq^{1-n},
 aq^{n+1} \cr};q,q \right] = \delta _{n,0}\;,\eqno(4.1.7)$$
where
$$\delta_{m,n} = \cases{1,&$m = n$,\cr
0,&$m\not= n$,\cr}$$
\noindent is the Kronecker delta function.  In \S 4.2 we will
need to use the $_4\phi_3$ summation formula
$$_4\phi_3\left[\matrix{a, qa^\half, -qa^\half, q^{-n}\cr
\noalign{\smallskip}
a^\half, -a^\half
, aq^{n+1} \cr};q,q^n\right] = \delta _{n,0}\;\eqno(4.1.8)$$
when $n = 0, 1, \ldots\,, $ which is the $b \to 0$ and the $b \to
\ty$ 
limit cases of (4.1.7).
The above derivation of (4.1.8) is substantially simpler than
that
in \S 2.3 of BHS, which used  the  
$q$-Pfaff--Saalsch\"utz formula (3.5.1) and the 
Bailey [1941] and Daum [1942] summation formula
$$_2\phi_1 (a,b;aq/b;q,-q/b) = {(-q;q)_\infty (aq,
aq^2/b^2;q^2)_\infty\over (aq/b,-q/b;q)_\infty}.\eqno(4.1.9)$$
Formula (4.1.9) is a $q$-analogue of Kummer's formula
$$_2F_1 (a,b;1+a-b;-1) = {\Gamma (1+a-b) \Gamma (1 + \half
a)\over \Gamma (1+a) \Gamma (1 + \half a-b)}.\eqno (4.1.10)$$

The series in (4.1.7) and (4.1.8) are the form
$$_{r+1}\phi_r\left[\matrix{a_1,
a_2,\ldots,a_{r+1}\cr 
b_1,\ldots,b_{r}\cr};q,z\right]\eqno(4.1.11)$$
in which the  the parameters satisfy the
relations
$$qa_1 = a_2b_1 = a_3b_2 = \cdots = a_{r+1}b_r .\eqno (4.1.12)$$
Such series are called {\it well-poised} series, and they are
called
{\it very-well-poised}\ if, in addition
$a_2 = qa_1^\half, a_3 = -qa_1^\half, $ which is the case for 
(4.1.7) and (4.1.8).
Since  very-well-poised series 
appear quite offen, to  simplify some of the displays
containing very-well-poised $_{r+1}\phi_r$ series, we will
frequently replace
$$_{r+1}\phi_r\left[\matrix{a_1, qa_1^\half, -qa_1^\half,
a_4,\ldots,a_{r+1}\cr a_1^\half, -a_1^\half,
qa_1/a_4,\ldots,qa_1/a_{r+1}\cr};q,z\right]\eqno(4.1.13)$$
\noindent with the  more compact notation
$$_{r+1}W_r \left(a_1; a_4, a_5, \ldots,
a_{r+1};q,z\right).\eqno(4.1.14)$$

Observe that the
expression on the right side  of (4.1.6) is {\it balanced}
 and {\it well-poised} since
$$\eqalignno{
&(ap)(bp)(cq)(aq/bc) = q (aq/b) (ap/c)(bcp)\cr
\noalign{\hbox{and}}
&(ap)q = (bp)
(aq/b) = (cq) (ap/c) = (aq/bc) (bcp).\cr}$$
Also, the part of the series on
the left side of (4.1.6) containing the $q$-shifted factorials is
{\it split-poised} in the sense that 
$aq = b(aq/b)$ and $c(ap/c) = (a/bc) (bcp) = ap.$
It is the first property and previous special cases, such as Wm.
Gosper's
$b\to 0$ limit case of
(4.1.6)
$$\sum^n_{k=0} {1-ap^kq^k\o 1-a} {(a;p)_k(c;q)_k\o (q;q)_k
(ap/c;p)_k} c^{-k}= {(ap;p)_n (cq;q)_n\o (q;q)_n (ap/c;p)_n} c^{-
n},\eqno (4.1.15)$$
that led the author to consider the sequence $s_k$ in (4.1.3)
and  obtain the factorization (4.1.4). 
See the discussion in Gasper [1989a, p. 259] . 

 This and the 
importance of the applications of (4.1.6) considered in Gasper
[1989a]
  led  Gasper and Rahman [1990b]
to generalize (4.1.3) to
$$s_k = {(ap,bp;p)_k (cq, ad^2q/bc;q)_k\o (dq, adq/b;q)_k (adp/c,
bcp/d;p)_k}\eqno (4.1.16)$$
for $k = 0, \pm 1, \pm 2, \ldots$, and observe  that 
$$\eqalignno{
\D t_k &= t_k - t_{k-1}\cr
&= {(ap, bp;p)_{k-1} (cq,ad^2q/bc;q)_{k-1}\o (dq, adq/b;q)_k
(adp/c, bcp/d;p)_k}\cr
&\quad\cdot \l\{ (1-ap^k) (1-bp^k) (1-cq^k) (1-ad^2q^k/bc)\r.\cr
&\qquad \l.- (1-dq^k)(1-adq^k/b) (1-adp^k/c) (1-bcp^k/d)\r\}\cr
&= {d(1-c/d) (1-ad/bc) (1-adp^kq^k) (1-bp^k/dq^k)\o (1-a) (1-b)
(1-c) (1-ad^2/bc)}\cr
&\quad\cdot {(a,b;p)_k (c,ad^2/bc;q)_kq^k\o (dq, adq/b;q)_k
(adp/c,
bcp/d;p)_k}&(4.1.17)\cr}$$
by (4.1.5), which, in view of (4.1.1),  gives the following
generalization
 of (4.1.6)
$$
\eqalignno{
&\sum^n_{k=-m} {(1-adp^kq^k)(1-bp^k/dq^k)\o (1-ad)(1-b/d)}
{(a,b;p)_k (c,ad^2/bc;q)_k\o (dq,adq/b;q)_k (adp/c, bcp/d;p)_k}
q^k\cr
&= {(1-a)(1-b)(1-c)(1-ad^2/bc)\o d(1-ad) (1-b/d) (1-c/d)(1-
ad/bc)}\cr
&\quad\cdot\l\{ {(ap, bp;p)_n (cq, ad^2q/bc;q)_n\o (dq,
adq/b;q)_n
(adp/c, bcp/d;p)_n} - {(c/ad, d/bc;p)_{m+1} (1/d, b/ad;q)_{m+1}\o
(1/c, bc/ad^2;q)_{m+1} (1/a, 1/b;p)_{m+1}}\r\},\cr
&&(4.1.18)\cr}$$
where we used $-m$ as the lower limit of summation instead of
$m.$  This has
the advantage that, by letting $m \to \ty$ in (4.1.18), we
immediately see that
if $|p| < 1$ and $|q| < 1,$ then
(4.1.18) tends to the  bibasic summation formula
$$\eqalignno{
&\sum^n _{k=-\ty} {(1-adp^kq^k) (1-bp^k/dq^k)\o (1-ad)(1-b/d)}
{(a,b;p)_k (c,ad^2/bc;q)_k\o (dq, adq/b;q)_k (adp/c, bcp/d;p)_k}
q^k\cr
&= {(1-a)(1-b)(1-c) (1-ad^2/bc)\o d(1-ad)(1-b/d)(1-c/d)(1-
ad/bc)}\cr
&\quad\cdot\l\{ {(ap, bp; p)_n (cq, ad^2q/bc;q)_n\o (dq,
adq/b;q)_n (adp/c, bcp/d;p)_n} - {(c/ad, d/bc;p)_\ty (1/d,
b/ad;q)_\ty\o (1/c, bc/ad^2;q)_\ty (1/a, 1/b;p)_\ty}\r\} \cr
&&(4.1.19)\cr}$$
for integer $n,$ and with $n$ replaced by $\ty$ . Some
generalizations
 of (4.1.18) are given in
Chu [1993] and Bhatnagar and Milne [1995].

Returning to (4.1.6), notice that when $c=q^{-n}$ it reduces to
$$\sum^n_{k=0} {(1-ap^kq^k)(1-bp^kq^{-k})\o (1-a)(1-b)}
{(a,b;p)_k (q^{-n}, aq^n/b;q)_k\o (q, aq/b;q)_k (apq^n, bpq^{-
n};p)_k} q^k = \delta_{n,0}\eqno (4.1.20)$$
when $n = 0, 1, \ldots \ $. 
This formula was derived independently by Bressoud [1988],
Gasper [1989a], and Krattenthaler [1995].
If we replace $n, a,b$ and $k$  by $n-m, ap^mq^m,
bp^mq^{-m}$ and $j-m$, respectively, (4.1.20) gives the rather
general orthogonality relation
$$\sum^n_{j=m}a_{nj}b_{jm} = \delta_{n,m}\eqno (4.1.21)$$
with
$$\eqalignno{
a_{nj} &= {(-1)^{n+j} (1-ap^jq^j)(1-bp^jq^{-j}) (apq^n, bpq^{-
n};p)_{n-1}\o (q;q)_{n-j} (apq^n, bpq^{-n};p)_j (bq^{1-
2n}/a;q)_{n-j}},\cr
&&\cr
b_{jm} &= {(ap^mq^m, bp^mq^{-m};p)_{j-m}\o (q, aq^{1+2m}/b;q)_{j-
m}} \l(- {a\o b} q^{1+2m}\r)^{j-m} q^{2{ j-m\choose
2}}.&\cr}$$

Hence, the triangular matrix $A = (a_{nj})$ is inverse
to the triangular matrix $B = (b_{jm}),$ where $j, m, n$ are
nonnegative integers. Since inverse matrices
commute, a calculation of the $jk^{\rm th}$ term of $BA$ leads to
the orthogonality relation
$$\eqalignno{
&\sum^{j-k}_{n=0} {(1-ap^kq^k)(1-bp^kq^{-k}) (ap^{k+1}q^{k+n},
bp^{k+1}q^{-k-n};p)_{j-k-1}\o (q;q)_n (q;q)_{j-k-n} 
(aq^{2k+n}/b;q)_{j-k-1}}\cr
&\cdot \l(1 - {a\o b} q^{2k+2n}\r) (-1)^n q^{n(j-k-1) + { j-k-
n \choose 2}} = \delta_{j,k} \ .&(4.1.22)\cr}$$
By replacing $j,n, a,b$ by $n+k, k, ap^{-k-1}q^{-k}, bp^{-
k-1}q^k$, respectively, we see that (4.1.22) is equivalent to
the bibasic summation formula
$$\l(1 - {a\o p}\r)\l(1 - {b\o p}\r) \sum^n_{k=0} {(aq^k, bq^{-
k};p)_{n-1} (1-aq^{2k}/b)\o (q;q)_k (q;q)_{n-k} (aq^k/b;q)_{n+1}}
(-1)^k q^{ k \choose 2} = \delta_{n,0}\eqno
(4.1.23)$$
when $n = 0, 1, \ldots \ $.
 Al-Salam and
Verma [1984] derived the $b\to 0$ limit case of (4.1.23)
$$\l(1 - {a\o p}\r) \sum^n_{k=0} {(aq^k;p)_{n-1} \o (q;q)_k
(q;q)_{n-k}}
(-1)^k q^{n- k \choose 2} = \delta_{n,0}\eqno
(4.1.24)$$
when $n = 0, 1, \ldots \ ,$
by using the fact that the $n^{\rm th}$
$q$-difference
of a polynomial in $q$ of degree less than $n$ is equal to zero.
For  applications  to $q$-analogues of Lagrange inversion,
general
 expansion formulas, and to the positivity of certain sums
(kernels), see
Gessel and Stanton [1986], Exercises 4.2 -- 4.4, and Gasper
[1989a, 1989b].
\bigskip
{\bf  4.2 Expansion, summation and transformation formulas.} 
Let $k$ and 
$n$ be  nonnegative
integers. In order to derive rather 
general summation and transformation formulas for $_{r+1}\phi_r$
series, it is efficient to first derive a rather general
expansion
formula that follows by using,  as in \S 2.2 of BHS,
the $q$-Pfaff--Saalsch\"utz formula (3.5.1)
 in the form 
$$\eqalignno{&_3\phi_2\left(q^{-k}, aq^k, aq/
bc; aq/b, aq/c;q,q\right)\cr
&= {(c, q^{1-k}/b;q)_k\over (aq/b,
cq^{-k}/a;q)_k}\;=\;{(b,c;q)_k\over (aq/b, aq/c;q)_k}\;
\left({aq\over bc}\right)^k & (4.2.1)\cr}$$ 
to obtain, for any sequence $\{v_k\}$,
$$\eqalignno{
&\som{k=0}{n} {(b,c,q^{-n};q)_k\over (q,
aq/b,aq/c;q)_k
}\;v_k\cr&= \som{k=0}{n} \som{j=0}{k} {(aq/bc, aq^k, q^{-k};q)_j
(q^{-n};q)_k\over (q, aq/b, aq/c;q)_j
(q;q)_k}\;q^j \left({bc\over aq}\right)^k v_k\cr
&= \som{j=0}{n} \som{i=0}{n-j} {\left(aq/bc, aq^{i+j},
q^{-i-j};q\right)_j \left(q^{-n};q\right)_{i+j}\over
(q, aq/b, aq/c;q)_j(q;q)_{i+j}}q^j \left({bc\over
aq}\right)^{i+j} v_{i+j}\cr
&= \som{j=0}{n} {\left(aq/bc, aq^j,
q^{-n};q\right)_j\over (q,aq/b, aq/c;q)_j}
\;(-1)^jq^{-{ j\mathstrut \choose 2}}\cr
&\quad\cdot
\som{i=0}{n-j}
{\left(q^{j-n}, aq^{2j};q\right)_i\over(q, a
q^j;q)_i} q^{-ij} \left({bc\over aq}\right) ^{i+j}v_{i+j}. &
(4.2.2)\cr}$$

 Setting
$$v_k = {(a,a_1,\ldots, a_r;q)_k\over
(b_1,b_2,\ldots, b_{r+1};q)_k}\;z^
k\eqno (4.2.3)$$
\noindent in (4.2.2) yields the desired expansion formula
$$\eqalignno{&_{r+4}\phi_{r+3}\left[\matrix{a, b, c, a_1,
a_2,\ldots, a_r, q^{-n} \cr
\noalign{\smallskip}
aq/b, aq/c ,b_1, b_2,\ldots, b_r,
b_{r+1}
\cr};q,z\right]\cr
&= \som{j=0}{n} {\left(aq/bc, a_1, a_2,\ldots,
a_r,q^{-n};q\right)_j\over(q, aq/b, aq/c, b_1,\ldots, b_r,
b_{r+1};q)_j}\;\left(-{bcz\over aq}\right)^jq^{-{ j\mathstrut
\choose 2}}(a;q)_{2j}\cr 
&\quad\cdot{}_{r+2}\phi_{r+1}\left[\matrix{
aq^{2j}, a_1q^j, a_2q^j, \ldots , a_rq^j,  q^{j-n}  \cr 
\noalign{\smallskip}
b_1q^j,
b_2q^j, \ldots, b_rq^j, b_{r+1}q^j \cr};q, {bcz\over
aq^{j+1}}\right].& (4.2.4)\cr}$$
This formula enables us to reduce the problem of
deriving a transformation formula for a $_{r+4}\phi_{r+3}$ series

in terms of a single series to that of
summing the 
$_{r+2}\phi_{r+1}$ series in (4.2.4) for some values of the
parameters.

By setting $a_1 = qa^\half,$ $a_2 =-qa^\half, b_1 = a^\half,$
$b_2 = -a^\half, b_{r+1} = aq^{n+1}$ and $a_k = b_k$, for
$ k = 3,4,\ldots, r$ in (4.2.4), we get
$$\eqalignno{
&_6\phi_5\left[\matrix{a, qa^\half, -qa^\half, b, c, q^{-n}\cr
\noalign{\smallskip}
 a^\half, -a^\half, aq/b, aq/c,  aq^{n+1}\cr};q,z\right]\cr
&= \som{j=0}{n} {(aq/bc, qa^\half,-qa^\half,
q^{-n};q)_j(a;q)_{2j}
\over (q, a^\half,-a^\half, aq/b, aq/c, aq^{n+1};q)_j}\;
\left(-{bcz\over aq}\right)^jq^{-{ j\mathstrut\choose
2}}\cr
&\quad\cdot{}_4\phi_3\left[\matrix{
aq^{2j}, q^{j+1}a^\half, -q^{j+1}a^\half, q^{j-n}\cr
\noalign{\smallskip}
 q^ja^\half, -q^ja^\half, a
q^{j+n+1}\cr};q, {bcz\o aq^{j+1}}\right].& (4.2.5)\cr}$$
If we now set $z=aq^{n+1}/bc$, then we can use (4.1.8)
to sum the above $_4\phi_3$ series; thus 
 deriving the summation formula 
$$\eqalignno{&_6\phi_5\left[\matrix{
a, qa^\half, -qa^\half, b, c, q^{-n}\cr
a^\half, -a^\half, aq/b, aq/c, aq^{n+1}\cr};q, {aq^{n+1}\o
bc}\right]\cr
&= {(aq/bc, qa^\half, -qa^\half, q^{-n};q)_n\;(a;q)_{2n}\over
(q,a^\half, -a^\half, aq/b, aq/c, aq^{n+1};q)_n}\; (-1)^n
q^{n(n+1)/2}\cr
&= {(aq,aq/bc;q)_n\over (aq/b,aq/c;q)_n},& (4.2.6)\cr}$$
which sums a terminating very-well-poised
$_6\phi_5$ series.

Similarly, from (4.2.4), we obtain
$$\eqalignno{
&_8\phi_7\left[\matrix{a, qa^{\half},
 -qa^{\half}, b, c, d, e, q^{-n}\cr
\noalign{\smallskip}
a^{\half}, -a^{\half}, aq/b, aq/c, aq/d, aq/e, aq^{n+1}\cr}
;q, {a^2q^{2+n}\o bcde}\right]\cr
&= \som {j=0}{n} {(aq/bc, qa^{\half}, -qa^{\half}, d, e, q^{-n}
;q)_j \; (a;q)_{2j}\over (q, a^{\half}, -a^{\half}, aq/b,
aq/c, aq/d, aq/e, aq^{n+1};q)_j}\; \left(- {aq^{n+1}\over
de}\right)^jq ^{-{ j\mathstrut\choose  2}}\cr
&\quad\cdot{} _6\phi_5\left[\matrix{
aq^{2j}, q^{j+1}a^{\half}, 
-q^{j+1}a^{\half}, dq^j, eq^j, q^{j-n}\cr
\noalign{\smallskip}
 q^ja^{\half}, -q^ja^{\half}, aq^{j+1}/d,
 aq^{j+1}/e, aq^{j+n+1}
\cr};q, {aq^{1+n-j}\o de}\right]& (4.2.7)\cr}$$
in which we can employ (4.2.6) to sum the $_6\phi_5$ series and
derive 
Watson's [1929] transformation
formula for a terminating very-well-poised $_8\phi_7$
series as a multiple of a terminating balanced  $_4\phi_3$
series:
$$\eqalignno{&_8\phi_7 \left[\matrix{
a,  qa^\half, -qa^\half, b, c, d, e, q^{-n}\cr
\noalign{\smallskip}
 a^{\half}, -a^{\half}, aq/b, aq/c, aq/d, aq/e,
 aq^{n+1}\cr};q, {a^2 q^{2+n}\o bcde}\right]\cr
&= {(aq, aq/de;q)_n\over (aq/d, aq/e;q)_n} \; _4\phi_3
\left[\matrix{q^{-n}, d, e, aq/bc\cr
 aq/b, aq/c, deq^{-n}/a\cr};q,q\right].&(4.2.8)\cr}$$

If $a^2q^{n+1} = bcde,$
the $_4\phi_3$ series in (4.2.8) becomes a terminating
balanced  $_3\phi_2$ series,
 which can be summed by the $q$-Pfaff--Saalsch\"utz formula to
derive 
 Jackson's [1921] summation formula for a 
terminating very-well-poised $_8\phi_7$ series
$$\eqalignno
{&_8\phi_7\left[\matrix{a, qa^{\half}, -qa^{\half}, b, c, d, e, 
q^{-n}\cr
\noalign{\smallskip}
a^{\half}, -a^{\half}, aq/b, aq/c, aq/d, aq/e, aq^{n+1}\cr}
;q,q\right]\cr
&= {(aq, aq/bc, aq/bd, aq/cd;q)_n\over (aq/b, aq/c, aq/d,
aq/bcd;q)_n}\, ,& (4.2.9)\cr}$$
where $a^2q^{n+1} = bcde.$
This formula is a $q$-analogue of 
Dougall's [1907] $_7F_6$ summation formula
$$\eqalignno{&_7F_6\left[\matrix{a,\ 1+\half a,\ b,\ c,\ d,\ e,\
-n \cr
\half a, 1+a-b, 1+a-c, 1+a-d, 1+a-e, 1+a+n \cr};1\right]\cr
&= {(1+a)_n (1+a-b-c)_n (1+a-b-d)_n (1+a-c-d)_n\over (1+a-b)_n
(1+a-c)_n (1+a-d)_n (1+a-b-c-d)_n},& (4.2.10)\cr}$$
where the series is 2-balanced, i.e, $1+2a+n = b+c+d+e.$
The reason that this series is  2-balanced instead of balanced is
 that the {\it appropriate} $q$-analogue of the term
$(1 + {1 \over 2}a)_k/({1\over 2}a)_k$
$= (a +2k)/a$ in the $_7F_6$ series is not $(qa^\half
;q)_k/(a^\half ;q)_k $ $= (1-a^\half q^k)/(1-a^\half )$ but
$(qa^{\half}, - qa^{\half};q)_k/(a^{\half}, -
a^{\half};q)_k =(1-aq^{2k})/(1-a)$, 
and  this introduces an additional $q$-factor in
the ratio of the products of the numerator and denominator
parameters. Krattenthaler [1995] removed some of the mystery in
the
factorization (4.1.5) by observing  that it is equivalent to the
$n=1$ case of Jackson's $_8\phi_7$ sum (4.2.9).  

Watson [1929] showed that the $b,c,d,e\rightarrow\infty $ limit
case
of his transformation formula (4.2.9)
and Jacobi's triple product identity can be used
to give a relatively
simple proof of the famous Rogers-Ramanujan identities:
$$\som{n=0}{\infty} {q^{n^{2}}\over (q;q)_n} = {(q^2,q^3,q^5;q^5)
_\infty\over (q;q)_\infty},\eqno (4.2.11)$$
$$\som{n=0}{\infty} {q^{n(n+1)}\over (q;q)_n} =
{(q,q^4,q^5;q^5)_\infty\over (q;q)_\infty},\eqno (4.2.12)$$
\noindent
where $|q| <1$. See Hardy [1940] for an early history of these
identities.

An important limit case of Jackson's summation formula (4.2.9) is
the
sum of a nonterminating very-well-poised $_6\phi_5$ series
$$\eqalignno
{&_6\phi_5\left[\matrix{a, qa^{\half}, -qa^{\half}, b, c, d\cr
\noalign{\smallskip}
a^{\half}, -a^{\half}, aq/b, aq/c, aq/d\cr}
;q, {aq\o bcd}\right]\cr
&= {(aq, aq/bc, aq/bd, aq/cd;q)_\infty\over (aq/b, aq/c, aq/d,
aq/bcd;q)_\infty}& (4.2.13)\cr}$$
with $|aq/bcd| < 1,$ which follows by letting
 $n\rightarrow \infty $ in (4.2.9).  
When $d=a^{\half}$  this formula reduces to 
$$\eqalignno{
&_4\phi_3 \left[\matrix{a, -qa^{\half}, b, c\cr
\noalign{\smallskip}
 -a^{\half}, aq/b, aq/c\cr};q, {qa^\half\o bc}\right]\cr
&= {(aq, aq/bc, qa^{\half}/b, qa^{\half}/c;q)_\infty\over
(aq/b, aq/c, qa^{\half}, qa^{\half}/bc;q)_\infty},&(4.2.14)
\cr}$$
\noindent
where $|qa^{\half}/bc|<1,$ which is a $q$-analogue of Dixon's 
[1903] formula for the sum of a well-poised $ _3F_2$ series
$$\eqalignno{&_3F_2\left[\matrix{a, b, c; 1 +a-b, 1+a-c
\cr};1\right]\cr
&= {\Gamma (1+{1\over 2}a) \Gamma (
1+a-b) \Gamma (1+a-c) \Gamma (1+ \half a-b-c)\over \Gamma (1+a)
\Gamma (1 + \half a-b) \Gamma (1+ \half a-c) \Gamma
(1+a-b-c)}, &(4.2.15)\cr}$$
where Re $(1+\half a-b-c) >0.$ 
\bigskip
{\bf  4.3 Bailey's transformation formulas and some integral 
representations.}  In this section we conclude the some of the
most 
important transformation formulas for very-well-poised series.
By rewriting Jackson's formula
(4.2.9)  in the form
$$\eqalignno{
&_8\phi_7
\left[\matrix{\lambda, q\lambda^\half ,-q\lambda^\half , \lambda
b/a, \lambda c/a, \lambda d/a, aq^m, q^{-m}\cr
\noalign{\smallskip}
\lambda^\half ,-\lambda^\half , aq/b, aq/c, aq/d, \lambda q^{1-
m}/a, \lambda q^{m+1}\cr};q,q\right]\cr
&= {(b,c,d,\lambda q;q)_m\over (aq/b, aq/c, aq/d, a/\lambda
;q)_m}& (4.3.1)\cr}$$
with $\lambda = qa^2/bcd$ 
and proceeding as in (4.2.2), one obtains a sum of a
terminating very-well-poised $_8\phi_7$ series that can be summed
with
Jackson's formula, yielding Bailey's [1929] 
transformation formula between  two 
terminating $_{10}W_9$ series
$$\eqalignno{&_{10}W_9 \left(a;b, c, d, e, f, 
\lambda aq^{n+1}/ef, q^{-n};q,q\right)\cr
&= {(aq, aq/ef, \lambda q/e, \lambda q/f;q)_n\over (aq/e, aq/f,
\lambda q/ef, \lambda q;q)_n}\cr
&\quad\cdot{}_{10}W_9 \left(\lambda;\lambda b/a, 
\lambda c/a, \lambda d/a, e, f,
\lambda aq^{n+1}/ef, q^{-n};q,q\right)& (4.3.2)\cr}$$
with $\lambda = qa^2/bcd,$ where for compactness
 we employed the $_{10}W_9$ notation 
for very-well-poised series.

 Watson's transformation formula
(4.2.8) follows from (4.3.2) by letting $b, c,$ or
$d\rightarrow \infty. $
 By taking the limit $n\rightarrow
\infty $ of (4.3.2) we obtain a transformation formula for
nonterminating
very-well-poised $_8\phi_7$ series
$$\eqalignno{&_{8}W_7 \left(a;b, c, d, e, f 
;q, \lambda q/ef \right)\cr
&= {(aq, aq/ef, \lambda q/e, \lambda q/f;q)_\infty \over (aq/e,
aq/f, \lambda q, \lambda q/ef;q)_\infty}\cr
&\quad\cdot{}_{8}W_7 \left(\lambda;\lambda b/a, 
\lambda c/a, \lambda d/a, e, f;q,aq/ef\right)& (4.3.3)\cr}$$
with $\lambda = qa^2/bcd,$
 where, for convergence,
$\max (|aq/ef |,\; |\lambda q/ef) <1.$
Bailey  iterated (4.3.2)  to obtain
$$\eqalignno{
&_{10}W_9\left(a;b,c,d,e,f,a^3q^{n+2}/bcdef, q^{-
n};q,q\right)\cr
&= {(aq, aq/de, aq/df, aq/ef;q)_n\over (aq/d, aq/e, aq/f,
aq/def;q)_n}\cr
&\quad\cdot{} _{10}W_9\left(def q^{-n-1}/a; aq/bc, d, e, f,
bdefq^{-n-
1}/a^2, cdefq^{-n-1}/a^2,q^{-n};q,q\right).\cr && (4.3.4)\cr}$$

It is clear that the $_{10}W_9$ on the left side of (4.3.4) tends
to the
$_8\phi_7$ series on the left side of (4.3.3) as $n\rightarrow
\infty $.  However, in trying to take the $n\rightarrow
\infty $ limit of the right side we run into the problem that the
terms near both ends of the series on the
right side of (4.3.4) are large compared to those in the middle
for large $n$, which prevents us from  directly taking the
term-by-term
limit.   Bailey overcame this difficulty by choosing $n$
to be an odd integer $2m+1,$  dividing the series on the
right into two halves, each containing $m+1$ terms, 
reversing the order of the second series, and then taking the
limit
as  $m\rightarrow
\infty $ to derive the transformation formula
$$\eqalignno{
&_8\phi_7\left[\matrix{a, qa^{\half}, -
qa^{\half}, b, c, d, e, f\cr
\noalign{\smallskip}
a^{\half}, -a^{\half}, aq/b, aq/c, aq/d, aq/e, aq/f\cr}
;q, {a^2q^2\o bcdef}\right]\cr
&= {(aq, aq/de, aq/df, aq/ef;q)_\infty\over (aq/d, aq/e, aq/f,
aq/def;q)_\infty}\;_4\phi_3\left[\matrix{aq/bc, d, e, f\cr
\noalign{\smallskip}
 aq/b, aq/c, def/a\cr};q,q\right]\cr
&\quad + {(aq, aq/bc, d, e, f, a^2q^2/bdef,\;
a^2q^2/cdef;q)_\infty\over
(aq/b, aq/c, aq/d, aq/e, aq/f, a^2q^2/bcdef, def/aq;q)_\infty}\cr
&\quad\cdot{} _4\phi_3\left[\matrix{aq/de, aq/df, aq/ef,
a^2q^2/bcdef\cr
\noalign{\smallskip}
a^2q^2/bdef, a^2q^2/cdef,
aq^2/def\cr};q,q\right],&(4.3.5)\cr}$$
provided $|a^2q^2/bcdef|<1$ when the $_8\phi_7$ series on the
left
side does
not terminate.  

 Al-Salam and Verma [1982] observed that (4.3.5) 
is equivalent to the $q$-integral formula
$$\eqalignno{
&\int_{a}^{b} {(qt/a, qt/b, ct,dt;q)_\infty\over (et, ft,
gt,ht;q)_\infty}\;d_qt\cr
&= b(1-q) {(q, bq/a, a/b, cd/eh, cd/fh, cd/gh, bc,
bd;q)_\infty\over (ae, af, ag, be, bf, bg, bh,
bcd/h;q)_\infty}\cr
&\cdot \;_8W_7 \left(bcd/hq; be, bf, bg, c/h,
d/h;q,ah\right),& (4.3.6)\cr}$$
\noindent
where $cd = abefgh$ and $|ah|<1$.
Setting $h=d$ in (4.3.6) 
and then replacing $g$ by $d$ gives Sears' [1951a] nonterminating
extension of the  $q$-Pfaff--Saalsch\"utz formula 
$$\eqalignno{&_3\phi_2\left[\matrix{a, b, c\cr
e, f\cr};q,q\right]\;
= {(q/e, f/a, f/b, f/c;q)\over (aq/e, bq/e, cq/e, f;q)_\infty}\cr
&- {(q/e, a, b, c, qf/e;q)_\infty\over (e/q, aq/e, bq/e, cq/e,
f;q)_\infty}\cr
&\cdot \;_3\phi_2\left[\matrix{aq/e, bq/e, cq/e\cr
\noalign{\smallskip}
q^2/e, qf/e\cr};q,q\right]& (4.3.7)\cr}$$
\noindent
with $ef = abcq$ in the equivalent
 $q$-integral form
$$\eqalignno{
&\int_{a}^{b} {(qt/a, qt/b,ct;q)_\infty\over (dt,
et,ft;q)_\infty}\;d_qt\cr
&= b(1-q) {(q, bq/a, a/b, c/d, c/e, c/f;q)_\infty\over (ad, ae,
af, bd, be, bf;q)_\infty},&(4.3.8)\cr}$$
where $c=abdef.$  

Rahman [1984] employed this $q$-integral to give
a rather simple proof of the 
Askey and Wilson [1985] $q$-beta integral formula
$$\eqalignno{&\int_{-1}^{1} {h(x; 1, -1, q^\half, -q^\half)\over
h(x; a, b, c, d)} {dx\over \sqrt{1-x^2}}\cr
&= {2\pi (abcd;q)_\ty\o (q,ab, ac, ad, bc, bd,
cd;q)_\ty} \ ,&(4.3.9)\cr}$$
\noindent
where
$$h(x;a_1,a_2,\ldots , a_m) \equiv  h(x;a_1,a_2,\ldots ,
a_m;q) = h(x;a_1) h(x;a_2)\cdots h(x;a_m),$$
$$h(x;a) \equiv h(x;a;q) = \prod_{n=0}^{\ty} (1-2axq^n +
a^2q^{2n})$$
\noindent
and
$$\max\l(|a|, |b|, |c|, |d|, |q|\r) < 1.$$
\noindent This is the integral that Askey and Wilson [1985] used
to derived
their orthogonality relation for the 
 polynomials
$$\eqalignno{p_n(x) &\equiv p_n(x;a,b,c,d|q)\cr
&= (ab, ac,
ad;q)_na^{-n}\;_4\phi_3\l[\matrix{q^{-n}, abcdq^{n-1}, 
ae^{i\ta}, ae^{-i\ta}\cr 
 ab, ac, ad \cr};q,q\r],&(4.3.10)\cr}$$
which are now called the {\it Askey-Wilson polynomials}.

In \S 2.11 of BHS the $q$-integral representation (4.3.6) is
applied
to derive  Bailey's [1936] 3-term
transformation formula for  $_8W_7$ series
$$\eqalignno{
&_8W_7(a,  b, c, d, e, f
;q, {a^2q^2/ bcdef})\cr
&= {(aq, aq/de, aq/df, aq/ef, eq/c, fq/c, b/a,
bef/a;q)_\infty\over (aq/d, aq/e, aq/f, aq/def, q/c, efq/c, be/a,
bf/a;q)_\infty}\cr
&\quad\cdot{}_8W_7(ef/c,  aq/bc, 
aq/cd, ef/a, e, f
;q, {bd/ a})\cr
&\quad + {b\over a} {(aq, bq/a, bq/c, bq/d, bq/e, bq/f, d, e,
f;q)_\infty\over
(aq/b, aq/c, aq/d, aq/e, aq/f, bd/a, be/a, bf/a,
def/a;q)_\infty}\cr
&\quad\cdot {(aq/bc, bdef/a^2, a^2q/bdef;q)_\infty\over (aq/def,
q/c, b^2q/a;q)_\infty}\cr
&\quad\cdot{}_8W_7(b^2/a,  b, bc/a, 
bd/a, be/a, bf/a;q,
{a^2q^2/ bcdef}),&(4.3.11)\cr}$$
\noindent
where $|bd/a|<1$ and $|a^2q^2/bcdef|<1,$ and it is pointed out
that 
the special case when $qa^2=bcdef$ gives
$$\eqalignno{
&_8W_7(a,  b, c, d, e, f;q,q)\cr
&\quad - {b\over a} {(aq, c, d, e, f, bq/a, bq/c, bq/d, bq/e,
bq/f;q)_\infty\over
(aq/b, aq/c, aq/d, aq/e, aq/f, bc/a, bd/a, be/a, bf/a,
b^2q/a;q)_\infty}\cr
&\quad\cdot{}_8W_7(b^2/a, 
b, bc/a, bd/a, be/a, bf/a;q,q)\cr
&= {(aq, b/a, aq/cd, aq/ce, aq/cf, aq/de, aq/df,
aq/ef;q)_\infty\over
(aq/c, aq/d, aq/e, aq/f, bc/a, bd/a, be/a,
bf/a;q)_\infty},& (4.3.12)\cr}$$
\noindent
where $qa^2 = bcdef$. This nonterminating extension of
Jackson's summation formula (4.2.9), which
is called Bailey's $_8\phi_7$ 
summation formula, 
can be written in the equivalent $q$-integral form
$$\eqalignno{
&\int_{a}^{b} {(qt/a, qt/b, t/a^\half , -t/a^\half ,
qt/c, qt/d, qt/e, qt/f;q)_\infty\over
(t, bt/a, qt/a^\half , -qt/a^\half , ct/a, dt/a, et/a,
ft/a;q)_\infty}\;d_qt\cr
&= {b(1-q) (q, a/b, bq/a, aq/cd, aq/ce, aq/cf, aq/de, aq/df,
aq/ef;q)_\infty\over (b, c, d, e, f, bc/a, bd/a, be/a,
bf/a;q)_\infty},& (4.3.13)\cr}$$
where $qa^2 = bcdef.$

Also, in \S 2.12 of BHS the $q$-integral  (4.3.6) was used to
give a short
derivation of  Bailey's [1947] 4-term
transformation formula for  $_{10}W_9$ series 
$$\eqalignno{
&_{10}W_9(a,  b, c, d, e, f, g, h ;q,q)\cr
&\quad + {(aq, b/a, c, d, e, f, g, h, bq/c, bq/d;q)_\infty\over
(b^2q/a, a/b, aq/c, aq/d, aq/e, aq/f, aq/g, aq/h, bc/a,
bd/a;q)_\infty}\cr
&\quad\cdot{(bq/e, bq/f, bq/g,
bq/h; q)_\infty\over
(be/a, bf/a, bg/a, bh/a;q)_\infty}\cr
&\quad\cdot{}_{10}W_9(b^2/a,  b, bc/a,  bd/a,  be/a, bf/a,  bg/a,
bh/a 
;q,q)\cr
&= {(aq, b/a, \lambda q/f, \lambda q/g, \lambda q/h, bf/\lambda,
bg/\lambda, bh/\lambda ;q)_\infty\over (\lambda q, b/\lambda,
aq/f, aq/g,
aq/h, bf/a, bg/a, bh/a;q)_\infty}\cr
&\quad\cdot{}_{10}W_9(\lambda,  b, \lambda c/a, \lambda d/a,
\lambda e/a, f, g, h ;q,q)\cr
&\quad + {(aq, b/a, f, g, h, bq/f, bq/g, bq/h, \lambda c/a,
\lambda
d/a;q)_\infty\over
(b^2q/\lambda, \lambda/b, aq/c, aq/d, aq/e, aq/f, aq/g, aq/h,
bc/a, bd/a;q)_\infty}\cr
&\quad \cdot {(\lambda e/a, abq/\lambda c, abq/\lambda d,
abq/\lambda
e;
q)_\infty\over
(be/a, bf/a, bg/a, bh/a ;q)_\infty}\cr
&\quad\cdot{} _{10}W_9(b^2/\lambda,  b, bc/a, bd/a, be/a,
bf/\lambda, bg/\lambda,
bh/\lambda;q,q) &(4.3.14)\cr}$$
and it was observed that (4.3.14)
 can be written in terms of the $q$-integrals in the 
compact form
$$\eqalignno{
&\int_{a}^{b} {(qt/a, qt/b, ta^{-\half}, -ta^{-\half}, qt/c,
qt/d,
qt/e, qt/f, qt/g, qt/h ;q)_\infty\over
(t, bt/a, qta^{-\half}, -qta^{-\half}, ct/a, dt/a, et/a, ft/a,
gt/a,
ht/a ;q)_\infty}\;d_qt\cr
&= {a\over\lambda} {(b/a, aq/b, \lambda c/a, \lambda d/a, \lambda
e/a, bf/\lambda, bg/\lambda, bh/\lambda ;q)_\infty\over
(b/\lambda, \lambda q/b, c, d, e, bf/a, bg/a, bh/a ;q)_\infty}\cr
&\quad\cdot{} \int_{\lambda}^{b} {(qt/\lambda,qt/b,
t\lambda^{-\half},
-t\lambda^{-\half}, aqt/c\lambda, aqt/d\lambda, aqt/e\lambda,
qt/f, qt/g, qt/h ;q)_\infty\over
(t, bt/\lambda, qt\lambda ^{-\half}, -qt\lambda^{-\half}, ct/a,
dt/a,
et/a, ft/\lambda, gt/\lambda, ht/\lambda
;q)_\infty}\;d_qt, \cr &&(4.3.15)\cr}$$
\noindent
where $\lambda = qa^2/cde$ and $a^3q^2 = bcdefgh$.
\bigskip
\centerline{\bf Exercises 4}
\bigskip
\item{4.1} Prove that
$$\eqalignno{
&\sum^n_{k=0} {(1-adp^kq^k) (1-bp^k/dq^k)\o (1-ad) (1-b/d)}
{(a,b;p)_k (q^{-n}, ad^2q^n/b;q)_k\o (dq, adq/b;q)_k (adpq^n,
bp/dq^n;p)_k} q^k\cr
&= {(1-d)(1-ad/b)(1-adq^n) (1-dq^n/b)\o (1-ad) (1-d/b) (1-dq^n)
(1-adq^n/b)}, \qquad n = 0,1,\ldots \;.\cr}$$
\medskip
\item{4.2} Use (4.1.24) to show that
 Euler's transformation
formula
$$\sum^\ty_{n=0} a_n b_nx^n = \sum^\ty_{k=0} (-1)^k {x^k\o k!}
f^{(k)} (x)\Delta^k a_0,$$
where
$$f(x) = b_0 + b_1x + b_2x^2+\cdots$$
and
$$\Delta^ka_0 = \sum^k_{j=0} (-1)^j {k\choose j} a_{k-j},$$
has the bibasic extension
$$\eqalignno{
\sum^\ty_{n=0} A_nB_n (xw)^n &= \sum^\ty_{k=0} (apq^k;p)_{k-1}
x^k \sum^k_{n=0} {(1-ap^nq^n) w^nA_n\o (q;q)_{k-n}
(apq^k;p)_n}\cr
&\cdot\sum^\ty_{j=0} {(ap^kq^k;p)_j\o (q;q)_j} B_{j+k} (-x)^j
q^{\ss j\mathstrut \choose\ss 2}.\cr}$$
See Al-Salam and Verma [1984].
\medskip
\item{4.3} Use (4.1.6) to derive the
Gasper [1989a] bibasic expansion formula
$$\eqalignno{
\sum^\ty_{n=0}A_nB_n{(xw)^n\o (q;q)_n} &= \sum^\ty_{n=0} {(1-\g
p^nq^n) (1-\sigma p^nq^n)\o (q;q)_n} (-x)^n q^{n+ {\ss n
\mathstrut\choose\ss 2}}\phantom{00000000} \cr
&\cr
&\quad\cdot \sum^\ty_{k=0} {1-\g\sigma^{-1}q^{2n+2k}\o (q;q)_k
(\g
pq^{n+k}, \sigma pq^{-n-k};p)_n} B_{n+k}x^k\cr
&\cr
&\quad\cdot \sum^n_{j=0}{(q^{-n};q)_j (\g\sigma^{-1} q^{n+j +
1};q)_{n+k - j - 1}\o (q;q)_j}\cr
&\cr
&\quad\cdot (\g pq^j, \sigma pq^{-j};p)_{n-1} A_j C_{j,n + k - j}
w^j
q^{n (j-n-k)},\phantom{00} \cr}$$
where $A_j, B_j, C_{j,k}$ are complex numbers such that the
series
converge absolutely and
$C_{j,0} = 1$, for $j = 0, 1, \ldots\ $.
\medskip
\item{4.4}  Show that if $p = q$, then the  $\sigma\to\ty$ limit
case of
Ex. 4.3 gives the expansion formula
$$\eqalignno{
&_{r+t}\phi_{s+u}\l[\matrix{a_R, c_T\cr b_S, d_U\cr} ;q,xw\r]\cr
&\cr
&= \sum^\ty_{j=0} {(c_T, e_K, \sigma, \g q^{j+1}/\sigma;q)_j\o
(q, d_U, f_M, \g q^j;q)_j} \l({x\o\sigma}\r)^j [(-1)^j
q^{\ss j\mathstrut \choose\ss 2}]^{u+m-t-k}\cr
\cr
&\quad\cdot\, _{t+k+4}\phi_{u+m+3}\l[\matrix{\g  q^{2j}/\sigma,
q^{j+1} \sqrt{\g/\sigma}, -q^{j+1}
\sqrt{\g/\sigma}, \sigma^{-1},\cr
\noalign{\smallskip}
q^j \sqrt{\g/\sigma}, -q^j\sqrt{\g/\sigma}, \g
q^{2j+1}, d_{\lower2pt\hbox{$\scriptstyle U$}}
 q^j,\cr}\r.\cr
&\qquad \l.\matrix{c_{\lower2pt\hbox{$\scriptstyle T$}}q^j, e_
{\lower2pt\hbox{$\scriptstyle K$}}q^j\cr
\noalign{\smallskip}
 f_{\lower2pt\hbox{$\scriptstyle M$}}q^j\cr};q,
xq^{j(u+m-t-k)}\r]\cr
&\cr
&\quad\cdot{} _{r+m+2}\phi_{s+k+2}\l[\matrix{q^{-j}, \g q^j, a_
{\lower2pt\hbox{$\scriptstyle R$}},
f_{\lower2pt\hbox{$\scriptstyle M$}}\cr \g q^{j+1}/\sigma,
q^{1-j}/\sigma, b_{\lower2pt\hbox{$\scriptstyle S$}}, e_
{\lower2pt\hbox{$\scriptstyle K$}}\cr} ;q,
wq\r],\cr}$$
where we employed the contracted notation of representing
$a_1,\ldots, a_r$ by \break $a_{\lower2pt\hbox{$\scriptstyle
R$}}$, etc.
\medskip
\item{4.5} Derive formulas (4.3.2) and (4.3.4).
\medskip
\item{4.6} Show that the  $q$-integral formulas (4.3.6) and
(4.3.8)
are equivalent to formulas (4.3.5) and (4.3.7), respectively.
\medskip
\item{4.7}  Use the $q$-integral representation (4.3.6)
to derive formulas (4.3.11) and (4.3.14).
\bigskip
\centerline{\bf References}
\bigskip \item{} 
Al-Salam, W.A. and Verma, A. [1982] {\it  Some remarks on
$q$-beta
integral,}  Proc. Amer. Math. Soc., {\bf 85}, 360--362.
\medskip \item{} 
Al-Salam, W.A. and Verma, A. [1984] {\it  On quadratic
transformations of basic series,}  SIAM J. Math. Anal.,
{\bf 15}, 414--420.
\medskip \item{} 
Andrews, G.E. [1969] {\it On a calculus of partition functions,}
Pacific J. Math., {\bf 31}, 555--562.
\medskip \item{} 
Andrews, G.E. [1986] $q$-{\it Series: Their
Development and Application in Analysis, Number Theory,
Combinatorics, Physics,
and Computer Algebra}, CBMS Regional Conference Lecture
Series, {\bf 66}, Amer. Math. Soc., Providence, R. I.
\medskip \item{}
Andrews, G.E. and Askey, R. [1977] {\it Enumeration of
partitions:
the role of Eulerian series and $q$-orthogonal polynomials,}
 Higher Combinatorics (M. Aigner, ed.), Reidel, Boston,
Mass., pp. 3--26. 
\medskip \item{}
Andrews, G.E. and Askey, R. [1978] {\it A simple proof of
Ramanujan's summation of the $_1\psi_1$},  Aequationes
Math., {\bf 18}, 333--337.
\medskip \item{} 
Askey, R. [1980] {\it  Ramanujan's extensions of the gamma and
beta
functions}, Amer. Math. Monthly, {\bf 87}, 346--359.
\medskip \item{}
Askey, R. and Ismail, M.E.H. [1983] {\it  A generalization of
ultraspherical polynomials,}  Studies in Pure Mathematics
(P. Erd\H os, ed.), Birkh\"auser, Boston, Mass., pp. 55--78. 
\medskip \item{}
Askey, R. and Roy, R. [1986] {\it More $ q$-beta integrals},
 Rocky Mtn. J. Math., {\bf 16}, 365--372.
\medskip \item{} 
Askey, R. and Wilson, J.A. [1985] {\it Some basic hypergeometric
polynomials that generalize Jacobi polynomials,}  Memoirs
Amer. Math. Soc., {\bf 319}.
\medskip \item{} 
Bailey, W.N. [1929] {\it  An identity involving Heine's basic
hypergeometric series,}  J. London Math. Soc., {\bf 4}, 254--257.
\medskip \item{} 
Bailey, W.N. [1935] {\it Generalized Hypergeometric Series},
Cambridge University Press, Cambridge, reprinted by
Stechert-Hafner, New York, 1964.
\medskip \item{}
Bailey, W.N. [1936] {\it Series of hypergeometric type which are
infinite in both directions,}  Quart. J. Math. (Oxford), {\bf 7},
105--115.
\medskip \item{}
Bailey, W.N. [1941] {\it A note on certain $q$-identities,}
Quart. J. Math. (Oxford), {\bf 12}, 173--175. 
\medskip \item{}
Bailey, W.N. [1947] {\it Well-poised basic hypergeometric
series,}
 Quart. J. Math. (Oxford), {\bf 18}, 157--166.
\medskip \item{}
Berndt, B.C. [1993] {\it Ramanujan's theory of theta-functions},
Theta Functions From the Classical to the Modern (M. Ram Murty,
ed.),
CRM Proceedings \& Lecture Notes, {\bf 1}, Amer. Math. Soc.,
Providence, R.I,
pp. 1--63.
\medskip \item{} 
Bhatnagar, G. and Milne, S.C. [1995] {\it Generalized bibasic 
hypergeometric series and their $U(n)$ extensions}, to appear.
\medskip \item{} 
Bressoud, D.M. [1988]  {\it The Bailey Lattice: an introduction,}

 Ramanujan Revisited (G. E. Andrews {\it et al}., eds.),
Academic Press, New York, pp. 57--67.
\medskip \item{} 
Cauchy, A.-L. [1843] {\it M\'emoire sur les fonctions dont
plusieurs
valeurs sont li\'ees entre elles par une \'equation lin\'eaire,
et sur diverses transformations de produits compos\'es d'un
nombre
ind\'efini de facteurs,}  C. R. Acad. Sci. Paris, T. XVII, p.
523,
{\it Oeuvres de Cauchy}, 1$^{\rm re}$
s\'erie, T. VIII, Gauthier-Villars, Paris, 1893, pp. 42--50.
\medskip \item{} 
Chu, W.C. [1993] {\it Inversion techniques and combinatorial
identities},
Bullettino U.M.I., {\bf 7}, 737--760.
 \medskip \item{}
Daum, J.A. [1942] {\it The basic analog of Kummer's theorem,}
Bull.
Amer. Math. Soc., {\bf 48}, 711--713. 
 \medskip \item{}
Dixon, A.C. [1903] {\it  Summation of a certain series,}  Proc.
London Math. Soc. (1), {\bf 35}, 285--289.
\medskip \item{}
Dougall, J. [1907] {\it On Vandermonde's theorem and some more
general expansions,}  Proc. Edin. Math. Soc., {\bf 25}, 114--132.
 \medskip \item{}
Fine, N.J. [1988] {\it Basic Hypergeometric Series and
Applications}, Mathematical Surveys and Monographs, Vol. 27,
Amer. Math. Soc.,
Providence, R. I. 
\medskip \item{} 
Gasper, G. [1987] {\it Solution to problem \#6497 ($
q$-Analogues of a gamma function identity, by R. Askey),}  Amer.
Math.
Monthly, {\bf 94}, 199--201.
\medskip \item{}
 Gasper, G. [1989a] {\it Summation, transformation, and expansion
formulas for bibasic series,}  Trans. Amer. Math.
Soc., {\bf 312}, 257--277.
\medskip \item{}
Gasper, G. [1989b]  $q$-{\it Extensions of Clausen's formula and
of the inequalities used by de Branges in his proof of the
Bieberbach, Robertson, and Millin conjectures,}  SIAM J. Math.
Anal., {\bf 20}, 1019--1034.
\medskip \item{}
Gasper, G. and Rahman, M. [1990a] {\it Basic Hypergeometric
Series}, Encyclopedia
of Mathematics and Its Applications, {\bf 35},
Cambridge University Press, Cambridge and New York.
\medskip \item{} 
Gasper, G. and Rahman, M. [1990b]  {\it An indefinite bibasic
summation formula and some quadratic, cubic, and quartic
summation and transformation formulas,}  Canad. J. Math., {\bf
42}, 1--27.
\medskip \item{} 
Gauss, C.F. [1813] {\it Disquisitiones generales circa seriem
infinitam ...,}  Comm. soc. reg. sci. G\"ott. rec., Vol.
II; reprinted in {\it Werke} {\bf 3} (1876), pp. 123--162.
\medskip \item{} 
Gessel, I. and Stanton, D. [1986]  {\it Another family of
$q$-Lagrange
inversion formulas,} Rocky Mtn. J. Math., {\bf 16}, 373--384.
\medskip \item{} 
Hahn, W. [1949] {\it \"Uber Polynome, die gleichzeitig zwei
verschiedenen Orthogonalsystemen angeh\"oren,}  Math.
Nachr., {\bf 2}, 263-278.
\medskip \item{} 
Hardy, G.H. [1940] {\it Ramanujan}, Cambridge University
Press, Cambridge; reprinted by Chelsea, New York, 1978.
\medskip \item{} 
Heine, E. [1846] {\it \"Uber die Reihe ...,} J. reine angew.
Math., {\bf 32}, 210--212.
\medskip \item{} 
Heine, E. [1847] {\it Untersuchungen \"uber die Reihe ...,}  J.
reine angew.  Math., {\bf 34}, 285--328.
\medskip \item{} 
Heine, E. [1878] {\it Handbuch der Kugelfunctionen, Theorie
und Anwendungen}, Vol. 1, Reimer, Berlin.
\medskip \item{} 
Ismail, M.E.H. [1977] {\it A simple proof of Ramanujan's
$_1\psi_1$ sum,} Proc. Amer. Math. Soc., {\bf 63}, 185--186.
\medskip \item{} 
Jackson, F.H. [1904] {\it A generalization of the functions
$\Gamma(n)$ and $ x^n$,}  Proc. Roy. Soc. London, {\bf 74},
64--72.
\medskip \item{} 
Jackson, F.H. [1910a] {\it Transformations of $ q$-series,}
Messenger of Math., {\bf 39}, 145--153.
\medskip \item{} 
Jackson, F.H. [1910b] {\it On $ q$-definite integrals,} 
Quart. J. Pure and Appl. Math., {\bf 41}, 193-203.
\medskip \item{} 
Jackson, F.H. [1921] {\it  Summation of $q$-hypergeometric
series,} Messenger of Math., {\bf 50}, 101--112.
\medskip \item{} 
Jackson, M. [1950] {\it On Lerch's transcendant and the basic
bilateral hypergeometric series $_2\psi_2$,}  J.
London Math. Soc., {\bf 25}, 189--196.
\medskip \item{} 
Jacobi, C.G.J. [1829] {\it Fundamenta Nova Theoriae Functionum
Ellipticarum,} Regiomonti. Sumptibus fratrum Borntr\"ager;
reprinted in {\it Gesammelte Werke} {\bf 1} (1881), 49--239,
Reimer,
Berlin.
\medskip \item{} 
Koornwinder, T.H. [1989] {\it Representations of the twisted
$ SU$(2) quantum group and some $ q$-hypergeometric
orthogonal polynomials,}  Proc. Kon. Nederl. Akad. Wetensch.
Series A,
{\bf 92}, 97--117.
\medskip \item{} 
Krattenthaler, C. [1995] {\it A new matrix inverse,} Proc. Amer.
Math. 
Soc., to appear.
\medskip \item{}
Pfaff, J.F. [1797] {\it Observationes analyticae ad L. Euler
Institutiones Calculi Integralis,} Vol. IV, Supplem. II et IV,
Historia de 1793,  Nova acta acad. sci. Petropolitanae,
{\bf 11} (1797), pp. 38--57.
\medskip \item{} 
Rahman, M. [1984] {\it  A simple evaluation of Askey and Wilson's
$q$-beta integral,}  Proc. Amer. Math. Soc., {\bf 92}, 413--417.
\medskip \item{} 
Ramanujan, S. [1915] {\it Some definite integrals,}  Messenger
of Math., {\bf 44}, 10--18.
\medskip \item{} 
Rogers, L.J. [1894] {\it Second memoir on the expansion of
certain
infinite products,}  Proc. London Math. Soc., {\bf 25}, 318--343.
\medskip \item{} 
Rogers, L.J. [1895] {\it Third memoir on the expansion of certain
infinite products,}  Proc. London Math. Soc., {\bf 26}, 15--32.
\medskip \item{}
Saalsch\"utz, L. [1890] {\it  Eine Summationsformel,} Zeitschr.
Math. Phys., {\bf 35}, 186--188.
\medskip \item{} 
Sears, D.B. [1951a] {\it Transformations of basic hypergeometric
functions of special type,}  Proc. London Math. Soc. (2), {\bf
52}, 467--483.
\medskip \item{} 
Sears, D.B. [1951b] {\it On the transformation theory of basic
hypergeometric functions,}  Proc. London Math. Soc. (2),
{\bf 53}, 158--180.
\medskip \item{}
Slater, L.J. [1966] {\it Generalized Hypergeometric Functions},
Cambridge University Press, Cambridge.
\medskip \item{} 
Thomae, J. [1869] {\it Beitr\"age zur Theorie der durch die
Heinesche
Reihe ...,}  J. reine angew. Math., {\bf 70}, 258--281.
\medskip \item{} 
Thomae, J. [1870] {\it Les s\'eries Hein\'eennes sup\'erieures,
ou
les s\'eries de la forme ...,}  Annali di Matematica Pura ed
Applicata,
{\bf 4}, 105--138.
\medskip \item{}
Venkatachaliengar, K. [1988] {\it Development of Elliptic
Functions According
to Ramanujan,} Tech. Rep. no. 2, Madurai Kamaraj University,
Madurai.
\medskip \item{} 
Watson, G.N. [1929] {\it  A new proof of the Rogers-Ramanujan
identities,}  J. London Math. Soc., {\bf 4}, 4--9.
\medskip \item{} 
Whittaker, E.T. and Watson, G.N. [1965] {\it A Course of
Modern Analysis}, 4th edition, Cambridge University Press,
Cambridge.

\bigskip
Department of Mathematics, Northwestern University, Evanston, IL
60208, USA

E-mail address: george@math.nwu.edu

\bye